\documentclass[11pt,a4paper]{article}
\title{A Hybrid High-Order Method for a Class of Strongly Nonlinear Elliptic Boundary Value Problems}

\author{Gouranga Mallik\footnote{Department of Mathematics, School of Advanced Sciences, Vellore Institute of Technology, Vellore 632014, India. Email. gouranga.mallik@vit.ac.inn}  and 
Thirupathi Gudi\footnote{Department of Mathematics, Indian Institute of Science, Bangalore 560012, India  Email. gourangam@iisc.ac.in}
}

\usepackage{amsmath,amsthm,amssymb,enumerate}
\usepackage{gensymb}
\usepackage[square,sort&compress,comma,numbers]{natbib}

\usepackage[sc]{mathpazo}
\usepackage{multirow} 

\usepackage{multicol}

\usepackage{newtxtext,newtxmath}
\usepackage{subfig}
\usepackage{graphicx}
\usepackage{epstopdf}
\usepackage{hyperref}
\usepackage{cancel}
\usepackage[top=30mm,bottom=30mm,foot=3mm,footskip=10mm,left=25mm,right=20mm]{geometry}
\usepackage{tikz}
\usepackage{caption}
\usetikzlibrary{shapes,calc}
\usepackage{verbatim}
\usepackage{mathrsfs}
\usepackage{algorithm}
\usepackage{accents}

\newtheorem{theorem}{Theorem}[section]
\newtheorem{lemma}[theorem]{Lemma}
\newtheorem{proof of lemma}[theorem]{Proof of Lemma}

\newtheorem{corollary}[theorem]{Corollary}

\theoremstyle{definition}

\newtheorem{remark}[theorem]{Remark}
\numberwithin{equation}{section}

\usepackage{amssymb,enumerate}
\usepackage{gensymb}
\usepackage[square,sort&compress,comma,numbers]{natbib}

\usepackage[sc]{mathpazo}
\usepackage{multirow} 

\usepackage{multicol}
\usepackage{newtxtext,newtxmath}
\usepackage{subfig}
\usepackage{graphicx}
\usepackage{epstopdf}
\usepackage{hyperref}
\usepackage{cancel}
\usepackage{tikz}
\usepackage{caption}
\usetikzlibrary{shapes,calc}
\usepackage{verbatim}
\usepackage{mathrsfs}
\usepackage{algorithm}
\usepackage{accents}
\newcommand{\dx}{{\rm\,dx}}

\newcommand{\ds}{{\rm\,ds}}
\newcommand{\dt}{{\rm\,dt}}

\newcommand{\fl}{\quad\forall}

\newcommand{\G}{1}

\newcommand{\mbR}{\mathbb R}

\newcommand{\mbP}{\mathbb P}

\newcommand{\cF}{\mathcal F}

\newcommand{\cT}{\mathcal T}
\newcommand{\cL}{\mathcal L}

\newcommand{\cN}{\mathcal N}

\newcommand{\fT}{\mathfrak{T}}

\newcommand{\ba}{\boldsymbol{a}}
\newcommand{\bz}{\boldsymbol{z}}
\newcommand{\bp}{\boldsymbol{p}}
\newcommand{\bq}{\boldsymbol{q}}
\newcommand{\bxi}{\boldsymbol{\xi}}

\newcommand{\tf}{\tilde{f}}
\newcommand{\tR}{\tilde{R}}

\newcommand{\lt}{L^2(\Omega)}

\newcommand{\Gtk}{\boldsymbol{G}_T^{k}}

\newcommand{\Utk}{\underline{U}_T^k}
\newcommand{\Ghk}{\boldsymbol{G}_h^{k}}

\newcommand{\Holders}{H\"{o}lder's }

\newcommand{\Gardings}{G$\mathring{\text{a}}$rding }
\newcommand{\GardingsType}{G$\mathring{\text{a}}$rding-type }

\newcommand{\CS}{Cauchy--Schwarz }
\newcommand{\etal}{\emph{et al.} }

\newcommand{\integ}{\int_\Omega}
\newcommand{\bfnTF}{\boldsymbol{n}_{TF}}
\newcommand{\bfn}{\boldsymbol{n}}
\newcommand{\sit}{\sum_{T\in\mathcal{T}_h}\int_T}

\newcommand{\Ctr}{C_{\text{tr}}}
\newcommand{\Ctrc}{C_{\text{tr},\text{c}}}

\begin{document}
	\date{\today}
	\maketitle
\begin{abstract}
In this article, we design and analyze a Hybrid High-Order (HHO) finite element approximation for a class of strongly nonlinear boundary value problems. We consider an HHO discretization for a suitable linearized problem and show its well-posedness using the \Gardings type inequality. The essential ingredients for the HHO approximation involve local reconstruction and high-order stabilization. We establish the existence of a unique solution for the HHO approximation using the Brouwer fixed point theorem and contraction principle. We derive an optimal order a priori error estimate in the discrete energy norm. Numerical experiments are performed to illustrate the convergence histories.
\end{abstract}
	
\noindent{\bf Key words:}
 Hybrid High-Order methods,  second-order nonlinear elliptic problems,  Brouwer fixed point theorem, error estimates.

\section{Introduction}\label{sec:intro}
There has been a growing interest in polytopal finite element methods of lower and higher-order polynomial approximations for partial differential equations. A non-exhaustive list includes the Hybridizable Discontinuous Galerkin method of \cite{Cock_Dong_09_HybridDG,Cock_Gopal_Laz_09_Unif,Piet_Ern_12_DG_book}, the Virtual Element method of \cite{Beira_Brezz_Cang_Man_Mar_13_VEM,Beira_Brezz_Mari_13_VEM_Elast,Brezz_Falk_Mari_14_VEM_Mixed}, the Weak Galerkin method of \cite{Mu_Wang_Ye_15_WeakGaler_Nonlin,Wang_Ye_13_Weak_Galerkin,Wang_Ye_14_Mixed_Weak_Galerkin}, the Gradient Discretization methods of \cite{Dron_Eym_Herb_GDM_16,Dron_Eym_GDM_18_book,Piet_Dron_Man_18_GradDis_Polytope}, the Multiscale Hybrid-Mixed method of \cite{Arya_Hard_Pare_Vale_13_Multiscale} and the Hybrid High-Order method of \cite{Piet_Ern_Lem_14_arb_local,Piet_Ern_15_mesh}. We refer to \cite{Piet_Jero_HHO_Book_20} for a thorough review of the literature on polytopal methods. The Hybrid High-Order (HHO) method has some specific features that distinguish it from the others. It is based on local polynomial reconstruction and complies with physics. The method is robust with respect to various physical parameters. The design is dimension-independent and suitable for local static condensation, which reduces the computational cost of the matrix solver.

The HHO method has some close connections with the Hybridizable Discontinuous Galerkin (HDG) method. It proposes a different stabilization than the HDG method to maintain the high-order convergence rate. The nonconforming Virtual Element Methods (ncVEM) choose the projection of virtual function in the stabilization, whereas the HHO method considers the reconstruction operator for the same. However, both methods achieve a similar rate of convergence. We refer to \cite{Cock_Piet_Ern_16_HHO_HDG} for detailed discussions on various relations of HDG and ncVEM with the HHO method.

HHO method in the lowest-order case falls in the family of the Hybrid Mixed Mimetic~\cite{Dron_Eym_Gall_Her_10_Mimetic}, which includes the Hybrid Finite Volume \cite{Eym_Gall_Herb_10_NonconfMesh}, the Mixed Finite Volume \cite{Dron_14_FV_Diff,Dron_Eym_06_MixedFV} and the Mixed-Hybrid Mimetic Finite Differences \cite{Brez_Lipn_Simo_05_Mim_Dif}. In~\cite{Lem_21_Bridging_HHO_VEM}, the author has bridged the HHO method with the virtual element method. We refer to \cite{Bone_Ern_14_poly,Dron_Eym_Gall_10_unified,Brezzi_Lip_Shas_05_Conv_Mimetic,Brezzi_Lip_Shas_07_polyhedral, Kuzn_Lip_Shas_04_Mimetic_poly} for related works.
We state some pivotal works on HHO methods for linear PDEs such as pure diffusion \cite{Piet_Ern_Lem_14_arb_local}, advection-diffusion \cite{Piet_Dron_Ern_15_HHO_Adv_Diff_Rea}, viscosity-dependent Stokes problem \cite{Piet_Ern_Link_16_HHO_Stokes} and interface problems \cite{Burman_Ern_18_HHO_Interface}, for nonlinear problems such as elliptic obstacle problem \cite{Cicu_Ern_Gudi_20_HHO_Obstacle},
a nonlinear elasticity with infinitesimal deformations \cite{Botti_Piet_Socha_17_HHO_Nonlin_Elast}, steady incompressible Navier Stokes equations \cite{Piet_Krell_18_NSE} and Leray-Lions operators \cite{Piet_Dron_16_Leray,Piet_Dron_Harn_21_Improved_Leray}.

In this article, we design and analyze HHO finite element approximation for the following class of strongly nonlinear partial differential equations (PDEs):
\begin{subequations}\label{intro_StrongNonlin_pde}
	\begin{align}
		-\nabla{\cdot}\ba(x,u,\nabla u)+f(x,u,\nabla u)&=0\quad\text{in}~~ \Omega,\\
		u(x)&=0\quad\text{on}~~ \partial\Omega,
	\end{align}
\end{subequations}
where $\Omega$ is a convex polytopal domain in $\mathbb{R}^d$, $d\in \{2,3\}$ with the Lipschitz boundary $\partial\Omega$. For the sake of simplicity, the homogeneous boundary condition is considered. We assume that $\ba(x,y,\bz):\bar{\Omega}\times \mbR\times\mbR^d\to\mbR^d$ and $f(x,y,\bz):\bar{\Omega}\times \mbR\times\mbR^d\to\mbR$ are twice continuously differentiable functions with all partial derivatives bounded and that \eqref{intro_StrongNonlin_pde} has a solution $u\in H^1_0(\Omega)$, see \cite{Bi_Gin_11_FV_strong,Chun_Chen_Yan_18_apost_twogrid}.  The linearized operator (namely, the Fr\'{e}chet derivative at $u$ in the direction $\psi$) is given by
\begin{equation}\label{lin_intro}
	\cL(u)\psi\equiv-\nabla{\cdot}\big{(}\ba_{\bz}(x,u,\nabla u)\nabla\psi+\ba_y(x,u,\nabla u)\psi\big{)}+f_{\bz}(x,u,\nabla u){\cdot}\nabla\psi + f_y(x,u,\nabla u)\psi,
\end{equation}
where $\ba_{\bz}=D_{\bz} \ba$ and $\ba_y=D_{y}\ba$ denote the derivatives of $\ba$ with respect to $\bz$ and $y$ respectively.
Following \cite{Xu_96_TwoGrid,Bi_Gin_11_FV_strong,Bi_Gin_14_Global_sup}, we assume the following two conditions:
\begin{enumerate}
	\item The matrix $\ba_{\bz}(x,u,\nabla u)$ is a symmetric and uniformly positive definite in $\bar{\Omega}$. That is, there exists a positive constant $\lambda_0$ such that $\lambda_0|\xi|^2\leq \xi^T \ba_{\bz}(x,u,\nabla u)\xi$ for $\xi\in\mbR^d$ and $x\in\bar{\Omega}$.
	\item The linearized operator $\cL(u): H^1_0(\Omega)\to H^{-1}(\Omega)$ is an isomorphism.
\end{enumerate}
This ensures that $u$ is an isolated solution to \eqref{intro_StrongNonlin_pde}. It can be observed that if $-\nabla{\cdot} \ba_y(x,u,\nabla u)+f_y(x,u,\nabla u)\geq 0$ then $\cL$ is an isomorphism (see \cite[Theorem~8.9]{Gil_Trud_book_01} and \cite{Xu_96_TwoGrid} for more details). 

Problems of the type \eqref{intro_StrongNonlin_pde} arise in several areas of applications, such as \cite{Bi_Gin_11_FV_strong,Gudi_NN_AKP_08_Strongly_nonlin}:
\begin{itemize}
	\item the equation of prescribed mean curvature
	\begin{equation*}
		\ba(x,u,\nabla u)=\left(1+|\nabla u|^2\right)^{-1/2}\nabla u,\text{ and } f(x,u,\nabla u)=f(x);
	\end{equation*}
	\item the subsonic flow of an irrotational, ideal, compressible gas
	\begin{equation*}
		\ba(x,u,\nabla u)=\left(1-\frac{\gamma-1}{2}|\nabla u|^2\right)^{1/(\gamma-1)}\nabla u\text{ with }\gamma>1 \text{ and }f(x,u,\nabla u)=f(x).
	\end{equation*}
\end{itemize}

We highlight some of the essential articles on finite element approximation for \eqref{intro_StrongNonlin_pde}. In \cite{Xu_96_TwoGrid}, Xu proved the existence of a unique finite element solution and derived optimal error estimates in the $L^p$- and $W^{1,p}$-norms under the assumption $H^1_0(\Omega)\cap W^{2,2+\epsilon}(\Omega)$ for some $\epsilon>0$. In \cite{Dem_06_loc_point_apost}, Demlow studied the residual-based pointwise a posteriori error estimates for finite element approximations. Gudi \etal \cite{Gudi_NN_AKP_08_Strongly_nonlin} and Bi \etal \cite{Chun_Chen_Yanp_15_strongNonlin_apost} studied the a priori and a posteriori error estimates for the $hp$-discontinuous Galerkin methods for \eqref{intro_StrongNonlin_pde}, respectively, under the assumption of $u\in H^1_0(\Omega)\cap H^{5/2}(\Omega)$ for $d=2$. We also refer to \cite{Bi_Gin_14_Global_sup,Bi_Gin_11_FV_strong,Chun_Chen_Yan_18_apost_twogrid} for various a priori and a posteriori error estimates for the problem. In \cite{Piet_Dron_16_Leray,Piet_Dron_Harn_21_Improved_Leray}, Di Pietro \etal designed and analyzed the HHO finite element approximation for the steady Leray--Lions
equation (where $f(x,u,\nabla u)=f(x)$) under the monotonicity and Lipschitz type of continuity assumptions on $\ba(x,u,\nabla u)$.


We briefly review some of the work on strongly nonlinear second-order PDEs. Gudi \etal  \cite{Gudi_AKP_07_DG_quasi,Gudi_NN_AKP_08_HpDG_Quasi} studied the existence and uniqueness of the discontinuous Galerkin (DG) and the local $hp$-DG finite element approximations for the following quasilinear problem of nonmonotone type:
\begin{align}\label{quasilin_pde_intro}
	-\nabla{\cdot}(a(x,u)\nabla u)&= f(x)\quad\text{in}~~ \Omega.
\end{align}
Bi \etal  \cite{Chun_Vict_13_Apost_DG_Quasi,Chun_Vict_09_APost_FV_Quasi,Chun_Ming_12_DGFV_Quasi,Song_Zhang_15_SupConv} studied various a priori and a posteriori error estimates for \eqref{quasilin_pde_intro}.
Recently, Gudi \etal \cite{TG_GM_TP_22_HHO_Quasi} analyzed the HHO finite element approximation for \eqref{quasilin_pde_intro} and proved the existence of a local unique discrete solution using the Brouwer fixed point theorem and the contraction principle.
Houston \etal \cite{Hous_Rob_Suli_05_Quasi} considered a one parameter family of $hp$-dG methods
for a class of quasilinear elliptic problems of the type:
\begin{align}\label{quasilin_mu}
	-\nabla{\cdot}(\mu(x,|\nabla u|)\nabla u)&= f(x)\quad\text{in}~~ \Omega,
\end{align}
where the coefficient function $\mu$ satisfies a monotone condition, see \cite{Hous_Rob_Suli_05_Quasi} for more details.


In this article, we analyze the HHO approximation for the strongly nonlinear problem \eqref{intro_StrongNonlin_pde} and establish an optimal order a priori error estimate in the discrete energy norm under the assumption $u\in H^1_0(\Omega)\cap H^3(\Omega)$. We use local reconstruction and high-order stabilization in the discrete formulation.  We establish the existence of a local unique discrete solution for the HHO approximation of \eqref{intro_StrongNonlin_pde}. We suitably define a nonlinear map and establish that the map possesses a ball to ball mapping and contraction properties. The fixed point of the non-linear map eventually is the solution to the discrete problem. As a consequence of the ball to ball and contraction properties, we obtain the error estimate in the energy norm. We follow some of the techniques of \cite{TG_GM_TP_22_HHO_Quasi}, where they consider $\ba(x,u,\nabla u)=a(x,u)\nabla u$ which leads to a linearized problem with scalar coefficient $a(x,u)$. In this article, the leading coefficient for the linearization \eqref{lin_intro} is a matrix $\ba_{\bz}(x,u,\nabla u)$, which depends on $u$ and $\nabla u$. This requires involved error analysis, and it possesses several additional difficulties.

The organization of the paper is as follows. Section~\ref{sec:intro} is introductory in nature. In Section~\ref{sec:HHO_Dis}, we introduce some notation and state some preliminary results related to HHO discretization.
In Section~\ref{sec:StrongNonLin}, we design and analyze the HHO approximation for the solution to the strongly nonlinear elliptic problem. In Section~\ref{sec:Numerics}, numerical experiments are performed to substantiate the theoretical results.

Throughout the paper, standard notation on Lebesgue and Sobolev spaces and their norms are employed.
For $K\subset \Omega$, the $L^2$-inner product on $L^2(K)$ is denoted by $(\bullet,\bullet)_{K}$ and $L^2$-norm by $\|\bullet\|_{K}$. We omit the subscript for the domain specification when $K=\Omega$. For the general $L^p$-space, we specify the appropriate domain and space in the definition of norm.
The standard seminorm and norm on $H^{s}(\Omega)$ (resp. $W^{s,p} (\Omega)$) for $s>0$ are denoted by $|\bullet|_{s}$ and $\|\bullet\|_{s}$ (resp. $|\bullet|_{s,p}$ and $\|\bullet\|_{s,p}$ ). The positive constants $C$ appearing in the inequalities denote generic constants, which do not depend on the meshsize. The notation $a\lesssim b$ means that there exists a generic constant $C$ independent of the meshsize such that $a \leq Cb$. We abbreviate $a\lesssim b\lesssim a$ by $a\approx b$.

\section{Hybrid High-Order discretization}\label{sec:HHO_Dis}

\subsection{Discrete setting}
Let $(\cT_h)_{h>0}$ be a sequence of refined meshes, where the parameter $h$ denotes the meshsize and goes to zero during the refinement process. For all $h>0$, we assume that the mesh $\cT_h$ covers $\Omega$ exactly and consists of a finite collection of non-empty disjoint open polyhedral cells $T$ such that $\overline{\Omega}=\cup_{T\in\cT_h} \overline{T}$ and $h=\max_{h\in\cT_h} h_T$, where $h_T$ is the diameter of $T$. A closed subset $F$ of $\Omega$ is defined to be a mesh face if it is a subset of an affine hyperplane $H_F$ with positive $(d-1)$-dimensional Hausdorff measure and if either of the following two statements holds true: (i) There exist $T_1(F)$ and $T_2(F)$ in $\cT_h$ such that $F\subset\partial T_1(F)\cap \partial T_2(F)\cap H_F$; in this case, the face $F$ is called an internal face; (ii) There exists $T(F)\in\cT_h$ such that $F \subset \partial T(F) \cap \partial\Omega \cap H_F$; in this case, the face $F$ is called a boundary face. The set of mesh faces is a partition of the mesh skeleton, that is, $\cup_{T\in\cT_h}\partial T=\cup_{F\in\cF_h}\bar{F}$, where $\cF_h:=\cF_h^i\cup \cF_h^b$ is the collection of all faces that is the union of the set of all internal faces $\cF_h^i$ and the set of all boundary faces $\cF_h^b$. Let $h_F$ denote the diameter of $F\in \cF_h$. For each $T\in \cT_h$, the set $F_T:=\{F\in \cF_h\, |\, F\subset\partial T\}$ denotes the collection of all faces contained in $\partial T$, $\bfn_T$ the unit outward normal to $T$ and we set $\bfnTF:=\bfn_T|_F$ for all $F\in \cF_h$. Following \cite[Definition 1]{Piet_Ern_15_mesh}, we assume that the mesh sequence $(\cT_h)_{h>0}$ is admissible in the sense that, for all $h>0$, $\cT_h$ admits a matching simplicial  submesh $\fT_h$ (i.e., every cell and face of $\fT_h$ is a subset of a cell and a face of $\cT_h$, respectively) so that the mesh sequence $(\cT_h)_{h>0}$ is shape-regular in the usual sense and all the cells and faces of $\cT_h$ have a uniformly comparable diameter to the cell and face of $\cT_h$ to which they belong. Owing to \cite[Lemma~1.42]{Piet_Ern_12_DG_book}, for $T\in\mathcal{T}_h$ and $F\in\mathcal{F}_T$, $h_F$ is comparable to $h_T$ in the sense that
\begin{align*}
	\varrho^{2}h_T\leq h_F\leq h_T,
\end{align*}
where $\varrho$ is the mesh regularity parameter. Moreover, there exists an integer $N_{\partial}$ depending on $\varrho$ and $d$ such that (see \cite[Lemma~1.41]{Piet_Ern_12_DG_book})
\begin{align*}
	\underset{T\in\mathcal{T}_h}{\text{max}}\, \text{card} (\mathcal{F}_T)\leq N_{\partial}.
\end{align*}

\noindent  Let $\mathbb{P}_d^l(T)$ be the polynomial space of degree at most $l$ on $T\in\cT_h$. There exist real numbers $\Ctr$ and $\Ctrc$ depending on $\varrho$ but independent of $h$ such that the following discrete and continuous trace inequalities hold for all $T\in\cT_h$ and $F\in\cF_T$ (see \cite[Lemma~1.46 and 1.49]{Piet_Ern_12_DG_book})
\begin{align}
	\|v\|_F&\leq \Ctr h_F^{-1/2} \|v\|_T\quad\forall v\in\mathbb{P}_d^l(T),\label{dis_trace}\\ 
	\|v\|_{\partial T}&\leq C_{\text{tr},\text{c}}(h_T^{-1}\|v\|_T^2+ h_T\|\nabla v\|_T^2)^{1/2}\quad\forall v\in H^1(T).\label{cts_trace_ineq}
\end{align}
Let $\pi_T^l$ be the $L^2$-orthogonal projector on $\mathbb{P}_d^l(T)$. There exists a real number $C_{\text{app}}$ depending on $\varrho$ and $l$ but independent of $h$ such that for all $T\in\cT_h$, the following holds (see \cite[Lemma~1.58 \& 1.59]{Piet_Ern_12_DG_book}): For all $ s\in \{1,\ldots,l+1\}$ and all $v\in H^s(T)$,
\begin{align}
	|v-\pi_T^l v|_{H^m(T)}+ h_T^{1/2}|v-\pi_T^l v|_{H^m(\partial T)}\leq C_{\text{app}} h_T^{s-m}|v|_{H^s(T)},\quad \forall m\in \{0, \ldots,s-1\},\label{proj_est}
\end{align}
where $|\bullet|_{H^m(\partial T)}$ denotes the facewise $H^m$-seminorm when the boundary $\partial T$ of an element $T\in\cT_h$ is written as a union of faces.

\subsection{Discrete spaces}
Let $k\geq 1$ be a fixed polynomial degree. Let $\mathbb{P}_d^k(T)$ be the space of polynomials of degree at most $k$ on the cell $T\in\cT_h$ and $\mathbb{P}_{d-1}^k(F)$ be the space of polynomial of degree at most $k$ on the face $F\in\cF_h$. For $T\in\cT_h$, the local space of degrees of freedom (DOFs) is defined by
\begin{align}\label{local_dofs}
	\underline{U}_T^k&:= \mathbb{P}_d^k(T)\times \left\{\underset{F\in\mathcal{F}_T}{\times} \mathbb{P}_{d-1}^k(F)\right\}.
\end{align}
The global space of DOFs is obtained by patching interface values in \eqref{local_dofs} as
\begin{align*}
	\underline{U}_h^k&:= \left\{\underset{T\in\mathcal{T}_h}{\times} \mathbb{P}_d^k(T)\right\}\times \left\{\underset{F\in\mathcal{F}_h}{\times} \mathbb{P}_{d-1}^k(F)\right\}.
\end{align*}
Imposing the zero boundary condition in the above discrete space $\underline{U}_h^k$, we define
\begin{align*}
	\underline{U}_{h,0}^k&:= \left\{\underline{v}_h=\left((v_T)_{T\in\cT_h},(v_F)_{F\in\cF_h}\right)\in\underline{U}_h^k  \,|\, v_F\equiv 0\fl F\in\mathcal{F}_h^b\right\}.
\end{align*}
Let $\pi_F^k$ be the $L^2$-orthogonal projector on $\mathbb{P}_{d-1}^k(F)$. Define a local interpolation operator $I_T^k:H^1(T)\to \underline{U}_T^k$ such that for all $v\in H^1(T)$,
	\begin{align}
		I_T^kv:= (\pi_T^k v, (\pi_F^kv)_{F\in\mathcal{F}_T}).\label{defn_interpolant}
	\end{align}
	The corresponding global interpolation operator $I_h^k:H^1(\Omega)\to \underline{U}_h^k$ is given by
	\begin{align*}
		I_h^kv:= ((\pi_T^kv)_{T\in\mathcal{T}_h}, (\pi_F^kv)_{F\in\mathcal{F}_h})\fl v\in H^1(\Omega).
\end{align*}
When applied to $H_0^1(\Omega)$, $I_h^k$ maps onto $\underline{U}_{h,0}^k$. 

We state a direct and reverse Lebesgue embedding result and refer to \cite[Lemma 5.1]{Piet_Dron_16_Leray} for proof. 
\begin{lemma}[Lebesgue embeddings]\label{lem_Sob_Inv_Ineq} Let $\cT_h$ be a regular mesh with $T\in\cT_h$. Let $k\in\mathbb{N}$ and $q,m\in [1,\infty]$. Then
	\begin{equation}
		\|v\|_{L^q(T)}\approx |T|^{\frac{1}{q}-\frac{1}{m}}\|v\|_{L^m(T)}\fl v\in \mathbb{P}_d^k(T).
	\end{equation}
\end{lemma}

The Sobolev exponent $p^*$ of $p$ is defined by
\begin{align*}
	p^*:=
	\begin{cases}
		\frac{dp}{d-p}\quad \text{ if } p<d,\\
		+\infty\quad \text{ if } p\geq d.
	\end{cases}
\end{align*}
We state a discrete Sobolev embedding from \cite[Proposition 5.4]{Piet_Dron_16_Leray} as follows. For $\underline{v}_h\in \underline{U}_h^k$, we understand $v_h\in L^2(\Omega)$ by $v_h|_T=v_T$. 
\begin{lemma}[Discrete Sobolev embeddings]\label{lem_dis_emb} Let $(\cT_h)_{h>0}$ be an admissible mesh sequence of $\Omega\subset \mbR^d$. Let $1\leq q\leq p^*$ if $1\leq p<d$ and $1\leq q<\infty$ if $p\geq d$. Then, there exists $C$ only depending on $\Omega, \varrho, k,q$ and $p$ such that
	\begin{equation*}
		\|v_h\|_{L^q(\Omega)}\leq C\|\underline{v}_h\|_{1,p,h}\fl\underline{v}_h\in \underline{U}_{h,0}^k, 
	\end{equation*}
	where $\displaystyle\|\underline{v}_h\|_{1,p,h}:=\left(\sum_{T\in\cT_h}\|\underline{v}_T\|_{1,p,T}^p\right)^{1/p}$ with
	\begin{equation}\label{dis_Lp_norm}
		\displaystyle \|\underline{v}_T\|_{1,p,T}:=\left(\|\nabla v_T\|^p_{L^p(T)^d}+\sum_{F\in\cF_T}h_F^{1-p}\|v_F-v_T\|_{L^p(T)}^p\right)^{1/p}.
	\end{equation} 
	In particular, 
	\begin{equation}
		\|v_h\|_{L^6(\Omega)}\leq C\|\underline{v}_h\|_{1,2,h}\fl\underline{v}_h\in \underline{U}_{h,0}^k.\label{dis_emb}    
	\end{equation}
\end{lemma}

\subsection{Local reconstructions and stabilization operators}
For $T\in\cT_h$, we define the local reconstruction operator $R_T^{k+1}:\underline{U}_T^k\to \mathbb{P}_d^{k+1}(T)$ such that for $\underline{v}_T=(v_T,(v_F)_{F\in\cF_T})$,
\begin{subequations}\label{recons_oper}
	\begin{align}
		(\nabla R_T^{k+1}\underline{v}_T,\nabla w)_T&=(\nabla v_T,\nabla w)_T+ \sum_{F\in\mathcal{F}_T}(v_F-v_T,\nabla w{\cdot}\bfnTF)_F,\label{hho3}\\
		\left(R_T^{k+1}\underline{v}_T,1\right)_T&= \left(v_T,1\right)_T,\label{hho9}
	\end{align}
\end{subequations}
where \eqref{hho3} is enforced for all $w\in \mathbb{P}_d^{k+1}(T)$.  
A global reconstruction operator $R_h^{k+1}:\underline{U}_h^k\to\mathbb{P}_d^{k+1}(\cT_h)$ is defined by $R_h^{k+1} \underline{v}_h|_T=R_T^{k+1} \underline{v}_T$.

We define a local gradient reconstruction $\Gtk:\Utk\to \mathbb{P}_d^k(T)^d$ such that for all $\underline{v}_T\in\underline{U}^k_T$,
\begin{align}\label{Grad_recons}
	(\Gtk\underline{v}_T,\boldsymbol{\tau})_T=(\nabla v_T,\boldsymbol{\tau})_T+\sum_{F\in\cF_T}(v_F-v_T,\boldsymbol{\tau}{\cdot}\bfnTF)_F\fl\boldsymbol{\tau}\in\mathbb{P}_d^k(T)^d.
\end{align}
Moreover, the following identity holds, see \cite[Lemma~4.10]{Piet_Jero_HHO_Book_20} for more details
\begin{align}\label{Grad_recons_proj}
	(\Gtk\underline{v}_T,\boldsymbol{\tau})_T=(\nabla v_T,\boldsymbol{\tau})_T+\sum_{F\in\cF_T}(v_F-v_T,(\pi_T^k\boldsymbol{\tau}){\cdot}\bfnTF)_F\fl\boldsymbol{\tau}\in L^1(T)^d.
\end{align}
The relation between $\Gtk$ and $R_T^{k+1}$ is established by taking $\boldsymbol{\tau}=\nabla w$ with $w\in\mbP_d^{k+1}(T)$ in \eqref{recons_oper} and comparing with \eqref{Grad_recons} as
\begin{align}
	(\Gtk\underline{v}_T-\nabla R_T^{k+1}\underline{v}_T,\nabla w)_T =0\fl w\in\mbP_d^{k+1}(T). 
\end{align}
In other words, $\nabla R_T^{k+1}\underline{v}_T$ is the $L^2$-orthogonal projection of $\Gtk\underline{v}_T$ on $\nabla \mbP_d^{k+1}(T)\subset  \mathbb{P}_d^k(T)^d$ and $\|\nabla R_T^{k+1}\underline{v}_T\|_T\leq \|\Gtk\underline{v}_T\|_T$.


The next lemma follows from \cite[Theorem~1.48]{Piet_Jero_HHO_Book_20} with  the trace inequality \eqref{cts_trace_ineq} and the approximation properties of an  elliptic projector $\pi_T^{1,k+1}$ since $R_T^{k+1}I_T^k v=\pi_T^{1,k+1}v$ for $v\in W^{1,1}(T)$.
\begin{lemma}[Approximation properties of $R_T^{k+1}I_T^k$]\label{lem_apprx_recons}
	There exists a real number $C>0$, depending on $\varrho$ but independent of $h_T$ such that for all $v\in H^{s+1}(T)$ for some $s\in\{0,1,\ldots,k+1\}$, 
		\begin{align}
			&\|v-R_T^{k+1}I_T^k v\|_T+ h_T^{1/2}\|v-R_T^{k+1}I_T^k v\|_{\partial T}+ h_T\|\nabla(v-R_T^{k+1}I_T^k v)\|_T \leq Ch_T^{s+1} |v|_{H^{s+1}(T)}.\label{eqn_recons_est}
		\end{align}
		For $s \in \{1,2,...,k+1\}$ and $v \in H^{s+1}(T)]$, we also have the approximation property
		\begin{align}\label{recon_approximation_1}
			h_T^{1/2}\|\nabla(v-R_T^{k+1}I_T^k v)\|_{\partial T}\leq Ch_T^{s} |v|_{H^{s+1}(T)}.
	\end{align}
\end{lemma}
The property $\Gtk I_T^k v=\pi_T^k(\nabla v)$ for $v\in W^{1,1}(T)$ and the approximation property for $L^2$ projector $\pi_T^k$ lead to
\begin{lemma}[Approximation properties of $\Gtk I_T^k$]\label{lem_apprx_Gtk}
	\cite[Lemma~3.24]{Piet_Jero_HHO_Book_20} There exists a real number $C>0$, depending on $\varrho$ but independent of $h_T$ such that for all $v\in H^{s+1}(T)$,
	\begin{align}
		\|\nabla v-\Gtk I_T^k v\|_{T}\leq Ch_T^{s} |v|_{H^{s+1}(T)} \quad\text{ for } s\in\{0,1,\ldots,k+1\}.\label{eqn_recons_Gtk}
	\end{align}
\end{lemma}

\section{Strongly nonlinear elliptic problem}\label{sec:StrongNonLin}
Let $\Omega$ be a bounded convex polytopal domain in $\mathbb{R}^d$, $d\in \{2,3\}$ with Lipschitz boundary $\partial\Omega$. In this article, we consider the HHO approximation for the strongly nonlinear elliptic boundary value problem:
\begin{subequations}\label{StrongNonlin_pde}
	\begin{align}
		-\nabla{\cdot}\ba(x,u,\nabla u)+f(x,u,\nabla u)&=0\quad\text{in}~~ \Omega,\\
		u&=0\quad\text{on}~~ \partial\Omega.
	\end{align}
\end{subequations}

For simplicity of notation, we often suppress $x$ in $a(x,u,\nabla u)$ and $f(x,u,\nabla u)$ when there is no confusion. Let $D:=\bar{\Omega}\times \mbR\times\mbR^d$.
We make the following assumptions on the problem~\eqref{StrongNonlin_pde}. 

\bigskip

\noindent\emph{{\bf Assumption N.1.} Nonlinear functions $f(x,y,\bz): D\to\mbR$ and $\ba(x,y,\bz):D\to\mbR^d$, are twice continuously differentiable with all their second-order derivatives bounded on $D$.}

\noindent\emph{{\bf Assumption N.2.} The derivative matrix $\left[a^{ij}(x,y,\bz)\right]_{i,j=1}^{d}=\left[\frac{\partial a_i}{\partial z_j}\right]_{i,j=1}^{d}$ for the coefficient function $\ba=(a_i)_{i=1}^d$ is symmetric.
There exist positive constants $\lambda_0$ and $\Lambda_0$ such that 
		\begin{equation}\label{cts_ellipticity}
			\lambda_0|\bxi|^2\leq \sum_{i,j=1}^d a^{ij}(x,u,\nabla u)\xi_i\xi_j\leq \Lambda_0|\bxi|^2 \fl x\in \Bar{\Omega} \text{ and } \xi\in \mbR^d.
		\end{equation} }
\emph{{\bf Assumption N.3.}  Assume that \eqref{StrongNonlin_pde} has a solution $u\in H^1_0(\Omega)$ with regularity $u\in H^3(\Omega)$.}

\begin{remark}
	For our subsequent error analysis, Assumption N.3 can be relaxed to $u\in H^1_0(\Omega) \cap H^{5/2}(\Omega)$ for $d=2$ and to $u\in H^1_0(\Omega)\cap H^{5/2+\epsilon}(\Omega),\,\epsilon>0$ for $d=3$. However, these require the approximation properties of \eqref{proj_est} and \eqref{eqn_recons_est} related to the projections $\pi_T^k$ and $\pi_F^k$ on fractional order Sobolev spaces, see \cite[Remark~1.49]{Piet_Jero_HHO_Book_20}. For simplicity of presentation, we kept our assumptions on integral Sobolev spaces.
\end{remark}

Using a suitable linearization, we design and analyze the HHO approximation for \eqref{StrongNonlin_pde}. The linearization of \eqref{StrongNonlin_pde} (namely, the Fr\'{e}chet derivative at $u$ in the direction $\psi$) is given by
\begin{align}\label{defn_cL}
	\cL(u)\psi\equiv-\nabla{\cdot}\big{(}\ba_{\bz}(u,\nabla u)\nabla\psi+\ba_y(u,\nabla u)\psi\big{)}+f_{\bz}(u,\nabla u){\cdot}\nabla\psi + f_y(u,\nabla u)\psi.
\end{align}

\noindent\emph{{\bf Assumption N.4.} The linearized operator $\cL(u): H^1_0(\Omega)\to H^{-1}(\Omega)$ is an isomorphism.}

\bigskip

In \cite{Bi_Gin_11_FV_strong}, the authors consider finite-volume-method for \eqref{StrongNonlin_pde} under the Assumption~N.1, N.2 and N.4, and establish optimal order a priori error estimates in the $W^{1,\infty}(\Omega)$ and $L^2$-norms under the regularity assumption $u\in W^{2,\infty}(\Omega)\cap H^3(\Omega)$. Gudi \etal \cite{Gudi_NN_AKP_08_Strongly_nonlin} and Bi \etal \cite{Chun_Chen_Yanp_15_strongNonlin_apost} derived the a priori and a posteriori error estimates for $hp$-discontinuous Galerkin methods for \eqref{intro_StrongNonlin_pde}, respectively, under the assumption of $u\in H^1_0(\Omega)\cap H^{5/2}(\Omega)$ for $d=2$.

If $-\nabla{\cdot}\ba_y(u,\nabla u)+f_y(u,\nabla u)\geq 0$ in addition to Assumptions~N.1 \& N.2, then the above Assumption~N.4 holds, see~\cite[Theorem~8.9]{Gil_Trud_book_01} and \cite{Xu_96_TwoGrid}.
Assumption~N.4 implies that the linearized problem: for given $\phi\in\lt$, find $\psi\in H^1_0(\Omega)$ such that
\begin{subequations}\label{lin_nonself}
	\begin{align}
		\cL(u)\psi&= \phi\quad\text{in}~~ \Omega,\\
		\psi&=0\quad\text{on}~~ \partial\Omega
	\end{align}
\end{subequations}
is well-posed. It can be observed that Assumption N.4 and an application of the open mapping theorem yield an a priori bound $\|\psi\|_{H^1(\Omega)}\lesssim \|\phi\|$, see \cite[Section 2.1]{Xu_96_TwoGrid}.
Since the domain $\Omega$ is convex,  the solution also satisfies the elliptic regularity $\|\psi\|_{H^2(\Omega)}\lesssim\|\phi\|$, see \cite[Lemma~2.1]{Xu_96_TwoGrid} and \cite{Grisvard_85_nonsmooth}. 
In the following sections, we consider an HHO approximation of the above linearized problem~\eqref{lin_nonself} and analyze the existence and uniqueness of the HHO approximation of \eqref{StrongNonlin_pde}.

\subsection{HHO approximations for a strongly nonlinear elliptic problem}
For $\underline{u}_h,\underline{v}_h\in\underline{U}_{h}^k$ define the discrete nonlinear form
\begin{align}
	\cN_h(\underline{u}_h; \underline{v}_h)&:=\sum_{T\in\mathcal{T}_h} \int_T
	\ba(u_T,\Gtk\underline{u}_T){\cdot}\Gtk\underline{v}_T\dx+ s_h(\underline{u}_h,\underline{v}_h) +\sit f(u_T,\Gtk\underline{u}_T)v_T\dx\label{defn_Bh},
\end{align}
where the above stabilization term $s_h(\underline{u}_h,\underline{v}_h)=\sum_{T\in\mathcal{T}_h}s_T(\underline{u}_T,\underline{v}_T)$ with the local contribution
	\begin{align}\label{stabilization}
		&s_T(\underline{u}_T,\underline{v}_T):=\frac{1}{h_T}\sum_{F \in \cF_T} \left(\pi_{F}^{k}(u_F-u_T-(R_T^{k+1}\underline{u}_T-\pi_T^{k}R_T^{k+1}\underline{u}_T)),\pi_{F}^{k}(v_F-v_T-(R_T^{k+1}\underline{v}_T-\pi_T^{k}R_T^{k+1}\underline{v}_T))\right)_F.
	\end{align} 
	We considered the scaling $h_T$ in place of $h_F$ for the above stabilization following the work of \cite{Jerome_Liam_21_hTscaling}.
The discrete HHO approximation of \eqref{StrongNonlin_pde} seeks $\underline{u}_h\in \underline{U}_{h,0}^k$ such that
\begin{align}\label{hho_dis_quasi}
	\cN_h(\underline{u}_h; \underline{v}_h)=0\quad \forall\underline{v}_h\in \underline{U}_{h,0}^k.
\end{align}
We establish the existence and uniqueness of a discrete solution to the above problem~\eqref{hho_dis_quasi} by a fixed point argument and the contraction result. We begin with a discrete linearized problem:
find $\underline{\psi}_h\in \underline{U}_{h,0}^k$ such that
\begin{align}\label{dis_lin_strong_nonlin}
	\cN^{\rm lin}_h(u;\underline{\psi}_h, \underline{v}_h)=(\phi,v_h)\quad \forall\underline{v}_h\in \underline{U}_{h,0}^k,
\end{align}
where we considered a linearization around the solution $u$ of \eqref{StrongNonlin_pde} and  for $\underline{\psi}_h, \underline{v}_h\in \underline{U}_{h,0}^k$,
\begin{align}
	\cN^{\rm lin}_h(u; \underline{\psi}_h, \underline{v}_h)&:= \sum_{T\in\mathcal{T}_h}\int_T \ba_{\bz}\Gtk\underline{\psi}_T{\cdot} \Gtk\underline{v}_T\dx+ s_h(\underline{\psi}_h,\underline{v}_h)+\sit\ba_y \psi_T{\cdot} \Gtk\underline{v}_T\dx\notag\\&\qquad+\sum_{T\in\mathcal{T}_h}\int_T f_{\bz}{\cdot}\Gtk\underline{\psi}_T v_T\dx+\sit f_y\psi_T v_T\dx.\label{defn_Bh_lin}
\end{align}
For the subsequent analysis, we also consider a fully discrete linearized form: for $\underline{w}_h,\underline{\psi}_h, \underline{v}_h\in \underline{U}_h^k$,
\begin{align}
	\tilde{\cN}^{\rm lin}_h(\underline{w}_h; \underline{\psi}_h, \underline{v}_h)&:= \sum_{T\in\mathcal{T}_h}\int_T \ba_{\bz}(w_T,\Gtk\underline{w}_T)\Gtk\underline{\psi}_T{\cdot} \Gtk\underline{v}_T\dx+\sit\ba_y(w_T,\Gtk\underline{w}_T)\psi_T{\cdot} \Gtk\underline{v}_T\dx+ s_h(\underline{\psi}_h,\underline{v}_h)\notag\\&\qquad+\sum_{T\in\mathcal{T}_h}\int_T f_{\bz}(w_T,\Gtk\underline{w}_T){\cdot}\Gtk\underline{\psi}_T v_T\dx+\sit f_y(w_T,\Gtk\underline{w}_T)\psi_T v_T\dx.\label{defn_Bh_lin_full}
\end{align}


\noindent Define a seminorm on $\underline{U}_h^k$ as follows:
\begin{align}\label{hho_gtk_norm}
	\|\underline{v}_h\|_{\G,h}^2:=\sum_{T\in\cT_h}\|\underline{v}_T\|_{\G,T}^2\text{ with } \|\underline{v}_T\|_{\G,T}^2:=\|\Gtk\underline{v}_T \|_T^2+\sum_{F\in\cF_T}\frac{1}{h_F}\|v_F-v_T\|_F^2.
\end{align}
Moreover, it is a norm in $\underline{U}_{h,0}^k$ owing to the zero boundary condition. It can be observed that the norm $\|\bullet\|_{1,2,h}$ in \eqref{dis_Lp_norm} is equivalent to $\|\bullet\|_{1,h}$  in $\underline{U}_h^k$.

In the next three lemmas, for simplicity of notation, we use $\ba_{\bz}, \ba_y, f_{\bz}$ and $f_y$ for $\ba_{\bz}(u,\nabla u), \ba_y(u,\nabla u), f_{\bz}(u,\nabla u)$ and $f_y(u,\nabla u)$ respectively, where there is no explicit role of $u$ and $\nabla u$. The following boundedness result can be obtained using the \CS inequality, the boundedness of $\ba_{\bz},\ba_y, f_{\bz},f_y$ and the definition of reconstructions $\Gtk,R_T^{k+1}$, see also \cite[Proposition 2.13]{Piet_Jero_HHO_Book_20}.
\begin{lemma}[Boundedness]
	For $\underline{u}_h, \underline{v}_h\in\underline{U}_h^k$, there exists a constant $C$ independent of meshsize $h$ such that 
	\begin{align}\label{bdd_ineq}
		\cN_h^{\rm lin}(u;\underline{u}_h,\underline{v}_h)\leq C \left(\|\underline{u}_h\|_{\G,h}+\| u_h\|\right)\left(\|\underline{v}_h\|_{\G,h}+\| v_h\|\right).
	\end{align}
\end{lemma}

We state and prove a \GardingsType inequality, which will be used to establish the existence of a solution to \eqref{dis_lin_strong_nonlin}.
\begin{lemma}[\GardingsType inequality]\label{lem:Garding}
	There exist two real numbers $C_1,C_2> 0$ independent of $h$ such that
	\begin{align}\label{garding_ineq}
		\cN_h^{\rm lin}(u;\underline{v}_h,\underline{v}_h)\geq C_1\|\underline{v}_h\|_{\G,h}^2- C_2\| v_h\|^2\quad \forall\underline{v}_h\in\underline{U}_{h}^k.
	\end{align}
\end{lemma}
\begin{proof}
	The first two terms of $\cN_h^{\rm lin}(u;\underline{v}_h,\underline{v}_h)$ in \eqref{defn_Bh_lin} are estimated using Assumption~N.2 and the lower bound of the stabilization of \cite[Proposition 2.13]{Piet_Jero_HHO_Book_20} as
	\begin{align}
		\sum_{T\in\mathcal{T}_h}\int_T \ba_{\bz}\Gtk\underline{v}_T{\cdot} \Gtk\underline{v}_T\dx+ s_h(\underline{v}_h,\underline{v}_h)\geq C \|\underline{v}_h\|_{\G,h}^2
	\end{align}
	for some positive constant $C$. The last three terms of $\cN_h^{\rm lin}(u;\underline{v}_h,\underline{v}_h)$
	are estimated using the \CS inequality as
	\begin{align*}
		&\sit\ba_y v_T{\cdot} \Gtk\underline{v}_T\dx\notag+\sum_{T\in\mathcal{T}_h}\int_T f_{\bz}{\cdot}\Gtk\underline{v}_T v_T\dx+\sit f_y v_T v_T\dx\leq \tilde{C}_1\| v_h\|\|\underline{v}_h\|_{\G,h}+\tilde{C}_2\| v_h\|^2
	\end{align*}
	for some positive constants $\tilde{C}_1,\tilde{C}_2$.
	The above two estimates lead to the required result
	\begin{align*}
		\cN_h^{\rm lin}(u;\underline{v}_h,\underline{v}_h)\geq C\|\underline{v}_h\|_{\G,h}^2-\tilde{C}_1\| v_h\|\|\underline{v}_h\|_{\G,h}-\tilde{C}_2\| v_h\|^2\geq C_1\|\underline{v}_h\|_{\G,h}^2- C_2\| v_h\|^2
	\end{align*}
	for some constants $C_1$ and $C_2$ independent of the meshsize $h$.
\end{proof}

In the following lemma, we prove the well-posedness of the linearized problem.  This is essential to propose a non-linear map, which is described in the next section.
\begin{lemma}\label{dis_lin_apriori}
	Adopt the aforementioned Assumptions~N.1--N.4. Assume $h$ is sufficiently small. For given $\xi\in L^2(\Omega)$, there exists a unique $\underline{\phi}_h\in \underline{U}_{h,0}^k$ such that
	\begin{align}\label{dis_dual}
		\cN_h^{\rm lin}(u;\underline{v}_h,\underline{\phi}_h)=(\xi,v_h)\fl \underline{v}_h\in\underline{U}_{h,0}^k.
	\end{align}
	Moreover, the solution $\underline{\phi}_h$ satisfies
	\begin{align}\label{dis_apriori}
		\|\underline{\phi}_h\|_{\G,h}\leq C\|\xi\|,
	\end{align}
	for sufficiently small $h$.
\end{lemma}
\begin{proof}
	First, we prove \eqref{dis_apriori}. Then the existence of a unique solution to (the finite dimensional system of equations) \eqref{dis_dual} follows immediately. The \Gardings type inequality \eqref{garding_ineq} with $\underline{v}_h=\underline{\phi}_h$ leads to 
	\begin{align*}
		C_1\|\underline{\phi}_h\|_{\G,h}^2&\leq \cN_h^{\rm lin}(u;\underline{\phi}_h,\underline{\phi}_h)+C_2\|\phi_h\|^2.
	\end{align*}
	Using \eqref{dis_dual} and the \CS inequality, we have
	\begin{align*}
		\cN_h^{\rm lin}( u;\underline{\phi}_h,\underline{\phi}_h)=(\xi,\phi_h)\leq \|\xi\|\|\phi_h\|\leq (\|\xi\|^2+\|\phi_h\|^2)/2.
	\end{align*}
	Combining the above two estimates, we obtain
	\begin{equation}\label{phi_gard}
		\|\underline{\phi}_h\|_{\G,h}\leq C_3\|\xi\|+C_4\|\phi_h\|. 
	\end{equation}
	We apply the Aubin-Nitche duality argument to estimate $\|\phi_h\|$. Consider the following auxiliary problem:
	\begin{subequations}\label{lin_aux_prob}
		\begin{align}
			-\nabla{\cdot}\left(\ba_{\bz}\nabla\psi+\ba_y\psi\right)+f_{\bz}{\cdot}\nabla\psi+f_y\psi&=\phi_h\quad\text{ in }\Omega,\\
			\psi&=0\quad\text{ in }\partial\Omega. 
		\end{align}
	\end{subequations}
	We recall the a priori bound for the solution $\psi\in H^1_0(\Omega)$ of \eqref{lin_aux_prob} from \eqref{defn_cL}--\eqref{lin_nonself}:
	\begin{equation}\label{H2_apriori_est}
		\|\psi\|_{H^2(\Omega)}\leq C\|\phi_h\|.
	\end{equation}
	Multiply \eqref{lin_aux_prob} by $\phi_h$ and integrate over $\Omega$ to obtain
	\begin{align}
		\|\phi_h\|^2&=-\integ \nabla{\cdot}(\ba_{\bz}\nabla \psi+\ba_y\psi)\phi_h\dx+\integ f_{\bz}{\cdot}\nabla\psi\phi_h\dx+\integ f_y\psi\phi_h\dx.\label{IBP_phi}
	\end{align}
	Since $\ba_{\bz}$ and $\ba_y$ are smooth and $\psi\in H^1_0(\Omega)\cap H^2(\Omega)$, we have the following two identities
	\begin{align}\label{IBP_bdd_zero}
		\sum_{T\in\cT_h}\sum_{F\in\cF_T}\int_F\phi_F\ba_{\bz}\nabla\psi{\cdot}\bfnTF\ds=0=\sum_{T\in\cT_h}\sum_{F\in\cF_T}\int_F\phi_F\ba_y\nabla\psi{\cdot}\bfnTF\ds,
	\end{align}
	see \cite[Corollary 1.19]{Piet_Jero_HHO_Book_20}.
	We apply the integration by parts on the first term of \eqref{IBP_phi} and use the identities \eqref{IBP_bdd_zero} and the definition of $\Gtk$ in \eqref{Grad_recons_proj} to obtain 
		\begin{align}
			&-\integ \nabla{\cdot}(\ba_{\bz}\nabla \psi+\ba_y\psi)\phi_h\dx=-\sit \nabla{\cdot}(\ba_{\bz}\nabla \psi+\ba_y\psi)\phi_h\dx\notag\\
			&=\sum_{T\in\cT_h}\left(\int_T \ba_
			{\bz}\nabla \psi{\cdot}\nabla \phi_T\dx+\sum_{F\in\cF_T}\int_F(\phi_F-\phi_T)\ba_{\bz}\nabla\psi{\cdot}\bfnTF\ds\right)\notag\\
			&\quad+\sum_{T\in\cT_h}\left(\int_T \ba_{y} \psi{\cdot}\nabla \phi_T\dx+\sum_{F\in\cF_T}\int_F(\phi_F-\phi_T)\ba_y\psi{\cdot}\bfnTF\ds\right)\notag\\
			&=\sum_{T\in\cT_h}\int_T \ba_{\bz}\nabla\psi{\cdot}\Gtk\underline{\phi}_T\dx+\sum_{T\in\cT_h}\sum_{F\in\cF_T}\int_F(\phi_F-\phi_T)\left(\ba_{\bz}\nabla\psi-\pi_T^k(\ba_
			{\bz}\nabla\psi)\right){\cdot}\bfnTF\ds\notag\\
			&\quad+\sum_{T\in\cT_h}\int_T \ba_y \psi{\cdot}\Gtk \underline{\phi}_T\dx+\sum_{T\in\cT_h}\sum_{F\in\cF_T}\int_F(\phi_F-\phi_T)\left(\ba_y\psi-\pi_T^k(\ba_y\psi)\right){\cdot}\bfnTF\ds.\notag\\
			&=\sum_{T\in\cT_h}\int_T \ba_{\bz}\Gtk I_T^k\psi{\cdot}\Gtk\underline{\phi}_T\dx+\sum_{T\in\cT_h}\int_T \ba_y\pi_T^k \psi{\cdot}\Gtk \underline{\phi}_T\dx\notag\\
			&\quad+\sum_{T\in\cT_h}\int_T \ba_{\bz}(\nabla\psi-\Gtk I_T^k\psi){\cdot}\Gtk\underline{\phi}_T\dx+\sum_{T\in\cT_h}\sum_{F\in\cF_T}\int_F(\phi_F-\phi_T)\left(\ba_{\bz}\nabla\psi-\pi_T^k(\ba_
			{\bz}\nabla\psi)\right){\cdot}\bfnTF\ds\notag\\
			&\quad+\sum_{T\in\cT_h}\int_T \ba_y (\psi-\pi_T^k\psi){\cdot}\Gtk \underline{\phi}_T\dx+\sum_{T\in\cT_h}\sum_{F\in\cF_T}\int_F(\phi_F-\phi_T)\left(\ba_y\psi-\pi_T^k(\ba_y\psi)\right){\cdot}\bfnTF\ds\notag\\
			&=:T_1+T_2+T_3+T_4+T_5+T_6.\label{term_az}
		\end{align} 
		The terms $T_3-T_6$ are estimated using the \CS inequality, the projection estimates of \eqref{proj_est} and Lemma~\ref{lem_apprx_Gtk} as
		\begin{align}
			T_3+T_4+T_5+T_6\leq Ch\|\psi\|_{H^2(\Omega)}\|\underline{\phi}\|_{1,h}.
	\end{align}
	The second and third terms of \eqref{IBP_phi} are controlled using the \CS inequality, the projection estimates of \eqref{proj_est} and Lemma~\ref{lem_apprx_Gtk} as follows
	\begin{align}
		& \integ f_{\bz}{\cdot}\nabla\psi\phi_h\dx+\integ f_y\psi\phi_h\dx\notag\\
		&=\integ f_{\bz}{\cdot}\Gtk (I_T^k\psi)\phi_h\dx+\integ f_{\bz}{\cdot}\left(\nabla\psi-\Gtk I_T^k\psi\right)\phi_h\dx\notag\\
		&\quad+\integ f_y\pi_h^k\psi\phi_h\dx +\integ f_y(\psi-\pi_h^k\psi)\phi_h\dx\notag\\
		&\leq \integ f_{\bz}{\cdot}\Gtk I_T^k\psi\phi_h\dx+\integ f_y\pi_h^k\psi\phi_h\dx+Ch\|\psi\|_{H^2(\Omega)}\|\underline{\phi}\|_{1,h}.\label{term_f}
	\end{align}
	Using the above estimates \eqref{term_az}--\eqref{term_f} in \eqref{IBP_phi}, we obtain
	\begin{align}
		\|\phi_h\|^2&\leq \cN_h^{\rm lin}(u;I_h^k\psi,\underline{\phi}_h)-s_h(I_h^k\psi,\underline{\phi}_h)+Ch\|\psi\|_{H^2(\Omega)}\|\underline{\phi}_h\|_{1,h}.\label{L2_est_phi}
	\end{align}
	Since  $s_h(I_h^k\psi,\underline{\phi}_h)\leq Ch\|\psi\|_{H^2(\Omega)}\|\underline{\phi}_h\|_{1,h}$ (see \cite[Equation~46]{Piet_Ern_Lem_14_arb_local} and 
	\begin{align}
		\cN_h^{\rm lin}(u;I_h^k\psi,\underline{\phi}_h)=\integ \xi\pi_h^k\psi\dx\leq \|\xi\|\|\pi_h^k\psi\|\leq\|\xi\|\|\psi\|_{H^2(\Omega)},
	\end{align}
	the above estimates and the a priori estimate \eqref{H2_apriori_est} in \eqref{L2_est_phi} lead to
	\begin{align}
		\|\phi_h\|\leq\|\xi\|+Ch\|\underline{\phi}_h\|_{1,h}.
	\end{align}
	This with \eqref{phi_gard} leads to $\|\underline{\phi}_h\|_{\G,h}\leq C\|\xi\|$ for sufficiently small $h$. This completes the proof.
\end{proof}


In the rest of the article, we use the following Taylor's formula in the integral form, see \cite{Gudi_NN_AKP_08_Strongly_nonlin,Bi_Cheng_Lin_16_Pointwise}: for $v\in\mbR$ and $\bp\in\mbR^d$ in terms of $u\in\mbR$ and $\bq\in\mbR^d$ 
	\begin{align}
		f(v,\bp)-f(u,\bq)&=f_y(u,\bq)(v-u)+f_{\bz}(u,\bq)(\bp-\bq)+\tilde{R}_f(v-u,\bp-\bq)\label{2nd_Order_Taylor}\\
		&=\tilde{f}_y(u,\bq)(v-u)+\tilde{f}_{\bz}(u,\bq)(\bp-\bq),\label{1st_order_Taylor}
	\end{align}
	where
	\begin{align*}
		\tf_y(u,\bq)=\int_0^1 f_y(u^t,\bq^t)\dt \text{ and } \tf_{\bz}(u,\bq)=\int_0^1 f_{\bz}(u^t,\bq^t)\dt.
	\end{align*}
	The remainder term $\tilde{R}_f$ in the above equation is given by, for $u^t=u+t(v-u), \bq^t=\bq+t(\bp-\bq)$,
	\begin{align}
		\tR_f(v-u,\bp-\bq)&=\tf_{yy}(u,\bq)(v-u)^2+2\tf_{y\bz}(u,\bq){\cdot}(\bp-\bq)(v-u)\notag\\
		&\qquad+(\bp-\bq)^T\tf_{\bz\bz}(u,\bq)(\bp-\bq),\label{eqn_Rf}
	\end{align}
	where
	\begin{align*}
		\tf_{yy}(u,\bq)&=\int_0^1(1-t)f_{yy}(u^t,\bq^t)\dt,\\
		\tf_{y\bz}(u,\bq)&=\int_0^1(1-t)f_{y\bz}(u^t,\bq^t)\dt\text{ and }\\
		\tf_{\bz\bz}(u,\bq)&=\int_0^1(1-t)f_{\bz\bz}(u^t,\bq^t)\dt.
	\end{align*}
	Similarly, the above Taylor's formula can be used for the function $\ba=(a_1,a_2)$ as:
	\begin{align}
		\ba(v,\bq)-\ba(u,\bq)&=\ba_y(u,\bq)(v-u)+\ba_{\bz}(u,\bq)(\bp-\bq)+\tR_{\ba}(v-u,\bp-\bq)\label{Res_Ra}\\
		&=\tilde{\ba}_{y}(u,\bq)(v-u)+\tilde{\ba}_{\bz}(u,\bq)(\bp-\bq),\notag
	\end{align}
	where
	\begin{align}
		\tR_{\ba}(v-u,\bp-\bq)=\left(\tR_{a_1}(v-u,\bp-\bq),\tR_{a_2}(v-u,\bp-\bq)\right)\label{eqn_Ra}
	\end{align}
	and
	\begin{align*}
		\tilde{\ba}_{y}(u,\bq)= \int_0^1\ba_y(u^t,\bq^t)\dt,\quad \tilde{\ba}_{\bz}(u,\bq)=\int_0^1\ba_{\bz}(u^t,\bq^t)\dt.
	\end{align*}
	Since $\ba$ and $f$ are twice continuously differentiable functions, all the above integral means involving second-order partial derivatives are bounded. That is, $\tilde{a}_{y},\tilde{a}_{\bz},\tilde{a}_{yy},\tilde{a}_{y\bz},\tilde{a}_{\bz y},\tilde{a}_{\bz\bz}$ and $\tilde{f}_{y},\tilde{f}_{\bz},\tilde{f}_{yy},\tilde{f}_{y\bz},\tilde{f}_{\bz y},\tilde{f}_{\bz\bz}\in L^\infty(D)$. Set 
	\begin{align}
		C_{\ba}:=\|\ba\|_{W^{2,\infty}(D )}, C_f:=\|f\|_{W^{2,\infty}(D)} \text{ and } C_{\ba,f}=\max\left\{C_{\ba},C_f\right\}.
	\end{align}


\subsection{fixed point formulation and contraction result}
In this section, we use fixed point arguments to establish the existence of a solution $\underline{u}_h\in \underline{U}_{h,0}^k$ of the above problem \eqref{hho_dis_quasi}. Local uniqueness is proved using the contraction principle. As a consequence of a fixed point result, an error estimate in the energy norm is deduced. Following the idea of \cite{GM_NN_CFEM,GM_NN_16_VKE_NCFEM,CC_GM_NN_19_VKE_DG}, we define a nonlinear map $\mu: \underline{U}_{h,0}^k\to \underline{U}_{h,0}^k$, which satisfies
\begin{align}
	\cN_h^{\rm lin}(u;I_h^ku-\mu(\underline{\theta}_h), \underline{v}_h)= \cN_h^{\rm lin}( u;I_h^ku-\underline{\theta}_h, \underline{v}_h)+ \cN_h(\underline{\theta}_h; \underline{v}_h)\fl \underline{v}_h\in \underline{U}_{h,0}^k.\label{nonlinear_map_mu}
\end{align}
The well-definedness of the map $\mu$ follows from the well-posedness of the linearized problem~\eqref{dis_lin_strong_nonlin}. We notice that any fixed point $\underline{\xi}_h$ (say) of $\mu$ satisfies the discrete problem~\eqref{hho_dis_quasi}. Now we proceed to prove the existence and uniqueness of a fixed point of the nonlinear map $\mu$. 
We make the following assumption throughout the section.

\noindent\textbf{Assumption N.5.} (Quasi-uniformity). We assume the admissible mesh sequence $(\cT_h)_{h>0}$ to be quasi-uniform, i.e., there exists a constant $C_Q$ independent of $h$ such that 
\begin{equation}\label{assump_quasi_uni}
	\max_{T\in\cT_h} h_T\leq C_Q\min_{T\in\cT_h}h_T.
\end{equation}


We propose some lemmas, which are used in the proof of the fixed point theorem.
\begin{lemma}\label{lem:diff_Lin}
	Let $u\in H^1_0(\Omega)\cap H^{r+2}(\cT_h)$ for $r\in\{0, 1,\ldots, k\}$. For $\underline{\theta}_h, \underline{v}_h\in \underline{U}_{h}^k$, it holds
	\begin{align}
		\left|\cN_h^{\rm lin}(u; \underline{\theta}_h, \underline{v}_h)-\tilde{\cN}_h^{\rm lin}(I_h^ku;\underline{\theta}_h, \underline{v}_h)\right|\leq Ch^{r+1-d/2}\|u\|_{H^{r+2}(\cT_h)}\|\underline{\theta}_h\|_{1,h}\|\underline{v}_h\|_{1,h}.\label{Bh_interpol}
	\end{align}
\end{lemma}
\begin{proof}
	From the definition of $\cN_h^{\rm lin}$ in \eqref{defn_Bh_lin} and $\tilde{\cN}_h^{\rm lin}$ in \eqref{defn_Bh_lin_full}, we have
	\begin{align*}
		&\cN_h^{\rm lin}(u; \underline{\theta}_h, \underline{v}_h)-\tilde{\cN}_h^{\rm lin}(I_h^ku;\underline{\theta}_h, \underline{v}_h)\notag\\
		&=\sum_{T\in\mathcal{T}_h}\int_T \left(\ba_{\bz}(u,\nabla u)-\ba_{\bz}(\pi_T^k u,\Gtk I_T^k u)\right)\Gtk\underline{\theta}_T{\cdot} \Gtk\underline{v}_T\dx+\sit(\ba_y(u,\nabla u)-\ba_y(\pi_T^k u,\Gtk I_T^k u)) \theta_T{\cdot} \Gtk\underline{v}_T\dx\notag\\&\qquad+\sum_{T\in\mathcal{T}_h}\int_T (f_{\bz}(u,\nabla u)-f_{\bz}(\pi_T^k u,\Gtk I_T^k u)){\cdot}\Gtk\underline{\theta}_T v_T\dx+\sit (f_y(u,\nabla u)-f_y(\pi_T^k u,\Gtk I_T^k u))\theta_T v_T\dx.
	\end{align*}
	The first term of the above equation is estimated by Taylor's formula~\eqref{1st_order_Taylor}, the generalized \Holders inequality, Lemma~\ref{lem_Sob_Inv_Ineq} and the definition of norm $\|\bullet\|_{1,h}$ in \eqref{hho_gtk_norm} as
	\begin{align*}
		&\sum_{T\in\mathcal{T}_h}\int_T \left(\ba_{\bz}(u,\nabla u)-\ba_{\bz}(\pi_T^k u,\Gtk I_T^k u)\right)\Gtk\underline{\theta}_T{\cdot} \Gtk\underline{v}_T\dx\\
		&\leq C_{\ba} \sum_{T\in\mathcal{T}_h}\left(\|u-\pi_T^ku\|_{T}+\|\nabla u-\Gtk I_T^k u\|_{T}\right)\|\Gtk\underline{\theta}_T\|_{L^4(T)}\|\Gtk\underline{v}_T\|_{L^4(T)}\\
		&\leq C C_{\ba} h^{r+1-d/2}\|u\|_{H^{r+2}(\cT_h)}\|\underline{\theta}_h\|_{1,h}\|\underline{v}_h\|_{1,h}.
	\end{align*}
	The remaining terms can be estimated in a similar way to obtain the desired result.
\end{proof}

The following three lemmas are essential to establish the fixed point result.
\begin{lemma}\label{lem:1st_remainder_est}
	Let $u\in H^1_0(\Omega)\cap H^{r+2}(\cT_h)$ for $r\in\{0, 1,\ldots, k\}$. For $\underline{v}_h\in \underline{U}_{h,0}^k$, the next three differences have the following estimates
	\begin{enumerate}[(i)]
		\item $\displaystyle\left|\sit \left(\ba(u,\nabla u)-\ba(\pi_T^ku,\Gtk I_T^k u)\right){\cdot}\nabla v_T\dx\right|\leq  CC_{\ba}h^{r+1}\|u\|_{H^{r+2}(\cT_h)}\|\underline{v}_h\|_{1,h}.$
		\item $\displaystyle \left|\sit \left(f(u,\nabla u)-f(\pi_T^ku,\Gtk I_T^ku)\right)v_h\dx\right|\leq  CC_{f}h^{r+1}\|u\|_{H^{r+2}(\cT_h)}\|\underline{v}_h\|_{1,h}.$
		\item $\displaystyle \left|\sum_{T\in\cT_h}\sum_{F\in\cF_T}\int_F( v_F- v_T)\left(\ba(u,\nabla u)-\pi_T^k\ba(u,\nabla u)\right){\cdot}\bfnTF\ds\right|
		\leq CC_{\ba}h^{r+1}\|\ba(u,\nabla u)\|_{[H^{r+1}(\cT_h)]^d}\|\underline{v}_h\|_{1,h}.$
	\end{enumerate}
\end{lemma}
\begin{proof}
The first two inequalities follow from the 1st-order Taylor’s formula~\eqref{1st_order_Taylor}, the \CS inequality, the property of the projection $\pi_T^k$ of \eqref{proj_est} and Lemma~\ref{lem_apprx_Gtk}. The third inequality follows from the \CS inequality and the estimation \eqref{proj_est} for the projection $\pi_T^k$.
\end{proof}

\begin{lemma}\label{lem:2nd_remainder_est}
For $\underline{\xi}_h,\underline{\chi}_h\in \underline{U}_{h}^k$ and $\underline{v}_h\in \underline{U}_{h,0}^k$, then we have the following bounds for the residuals:
	\begin{align*}
		&\left|\sit \tilde{R}_{f}(\xi_T,\Gtk\underline{\chi}_T)\,v_T\dx\right|\leq CC_f\big{(} \|\xi_h\|_{L^{4}(\Omega)}^2+\|\xi_h\|_{L^{4}(\Omega)}\|\Gtk\underline{\chi}_h\|_{L^{4}(\Omega)}+\|\Gtk\underline{\chi}_h\|_{L^{4}(\Omega)}^2\big{)}\|v_h\|.  
	\end{align*}
	and 
	\begin{align*}
		\left|\sum_{T\in\cT_h} \int_T \tilde{R}_{\ba}(\xi_T,\Gtk\underline{\chi}_T){\cdot}\Gtk\underline{v}_T\dx \right|\leq C C_{\ba}\big{(} \|\xi_h\|_{L^{4}(\Omega)}^2+\|\xi_h\|_{L^{4}(\Omega)}\|\Gtk\underline{\chi}_h\|_{L^{4}(\Omega)}+\|\Gtk\underline{\chi}_h\|_{L^{4}(\Omega)}^2\big{)}\|\Gtk\underline{v}_h\|.
\end{align*}
\end{lemma}
\begin{proof}
The proof follows from the definition of $\tilde{R}_{f}$ and $\tilde{R}_{\ba}$ with the generalized \Holders inequality.
\end{proof}

The next result follows from Lemma~\ref{lem:2nd_remainder_est} and Lemma~\ref{lem_Sob_Inv_Ineq}.
\begin{corollary}\label{cor:2nd_remainder_est}
For $\underline{\xi}_h\in \underline{U}_{h}^k$ and $\underline{v}_h\in \underline{U}_{h,0}^k$, the following bounds hold:
\begin{align*}
	\left|\sit \tilde{R}_{f}(\xi_T,\Gtk\underline{\xi}_T)\,v_T\dx\right|\leq C_f (\max_{T\in\cT_h} h_T^{-d/2})\|\underline{\xi}_h\|_{1,h}^2\|v_h\|
\end{align*}
and 
\begin{align*}
	\left|\sum_{T\in\cT_h} \int_T \tilde{R}_{\ba}(\xi_T,\Gtk\underline{\xi}_T){\cdot}\Gtk\underline{v}_T\dx \right|&\leq C_{\ba} (\max_{T\in\cT_h} h_T^{-d/2}) \|\underline{\xi}_h\|_{1,h}^2\|\Gtk\underline{v}_h\|.
\end{align*}
\end{corollary}

\begin{lemma}\label{lem:lin_est}
The following estimate for the linearization holds true
\begin{align}
	\cN_h(\underline{\theta}_h;\underline{v}_h)-\cN_h(I_h^k u;\underline{v}_h)+\tilde{\cN}_h^{\rm lin}(I_h^ku;I_h^ku-\underline{\theta}_h, \underline{v}_h)\leq C C_{\ba,f} h^{-d/2}\|\underline{\theta}_h-I_h^k u\|_{1,h}^2\|\underline{v}_h\|_{1,h}.
\end{align}
\end{lemma}

\begin{proof}
The definitions of $\cN_h$ and $\tilde{\cN}_h^{\rm lin}$ and Taylor's formula~\eqref{2nd_Order_Taylor} lead to the following identity
\begin{align}
	&\cN_h(\underline{\theta}_h;\underline{v}_h)-\cN_h(I_h^k u;\underline{v}_h)
	=\sum_{T\in\mathcal{T}_h} \int_T
	\ba(\theta_T,\Gtk\underline{\theta}_T){\cdot}\Gtk\underline{v}_T\dx-\sum_{T\in\mathcal{T}_h} \int_T
	\ba(\pi_T^k u,\Gtk I_T^k u){\cdot}\Gtk\underline{v}_T\dx\notag\\
	&\qquad +\sit f(\theta_T,\Gtk\underline{\theta}_T)v_T\dx-\sit f(\pi_T^k u,\Gtk I_T^k u)v_T\dx+s_h(\underline{\theta}_h-I_h^k u,\underline{v}_h)\notag\\
	&=\tilde{\cN}_h^{\rm lin}(I_h^ku;\underline{\theta}_h-I_h^ku, \underline{v}_h)+\sit \tilde{R}_{\ba}(\theta_T-\pi_T^k u,\Gtk(\underline{\theta}_T-I_T^k u)){\cdot}\Gtk\underline{v}_T\dx\notag\\
	&\quad+\sit \tilde{R}_{f}(\theta_T-\pi_T^k u,\Gtk(\underline{\theta}_T-I_T^k u))\,v_T\dx.
\end{align}
Estimates of Corollary~\ref{cor:2nd_remainder_est} and Assumption N.5 lead to the required result.
\end{proof}

Define a ball of radius $R$ with center at $I_h^k u$ as
\begin{equation*}
B(I_h^ku;R):=\left\{\underline{\theta}_h\in\underline{U}_{h,0}^k \text{ such that } \|I_h^k u-\underline{\theta}_h\|_{1,h}\leq R\right\},
\end{equation*}
and recall Assumption N.1--N.5 for the following result.
\begin{theorem}[fixed point result]\label{thm_fixed_point}
Let $u\in H^1_0(\Omega)$ be a solution to \eqref{StrongNonlin_pde}. Assume $u\in H^{r+2}(\cT_h)$ and $\ba(x,y,\bz)$ to be $(r+1)$ times continuously differentiable with respect to $x$, for some $r\in\{1,\ldots,k\}$.  Adopt the aforementioned Assumptions~N.1--N.5. For a sufficiently small meshsize $h$, there exists positive $R(h)$ such that the nonlinear map $\mu: \underline{U}_{h,0}^k\to \underline{U}_{h,0}^k$ defined in \eqref{nonlinear_map_mu} maps from the ball $B(I_h^ku;R(h))$ to itself. Moreover, $\mu$ has a fixed point in $B(I_h^ku;R(h))$  with a radius $R(h):=\tilde{C}h^{r+1}$ for some positive constant $\tilde{C}$ independent of the meshsize.
\end{theorem}

\begin{proof}
	From Lemma~\ref{lem:Garding}, we have 
	\begin{align}\label{Garding_infsup_0}
		C_1\|\underline{w}_h\|_{1,h}^2\leq \cN_h^{\rm lin}(u;\underline{w}_h, \underline{w}_h)+ C_2\|w_h\|^2\fl \underline{w}_h\in \underline{U}_{h,0}^k.
	\end{align}
	Using the inequality $\| w_h\|\leq C \|\underline{w}_h\|_{1,h}$ for $\underline{w}_h\in \underline{U}_{h,0}$ obtained from  Lemma~\ref{lem_dis_emb} and the \GardingsType inequality \eqref{Garding_infsup_0}, we have
	\begin{align}\label{Garding_infsup}
		C_1\|\underline{w}_h\|_{1,h}\leq \sup_{\underline{v}_h\in \underline{U}_{h,0}^k,\, \|\underline{v}_h\|_{1,h}=1}\cN_h^{\rm lin}(u;\underline{w}_h, \underline{v}_h)+ \tilde{C}_2\|w_h\|\fl \underline{w}_h\in \underline{U}_{h,0}^k,
	\end{align}
	for some positive constant $\tilde{C}_2$. Choose $\underline{w}_h=I_h^ku-\mu(\underline{\theta}_h)$ in the above equation. We understand $\|\underline{w}_{h}\|_{L^2}$ by $\|w_{h}\|$. Then, there exists $\underline{v}_h$ with $\|\underline{v}_h\|_{1,h}=1$ such that
\begin{align*}
	&C_1\|I_h^ku-\mu(\underline{\theta}_h)\|_{1,h}\leq \cN_h^{\rm lin}(u; I_h^ku-\mu(\underline{\theta}_h), \underline{v}_h)+ \tilde{C}_2\|I_h^ku-\mu(\underline{\theta}_h)\|_{L^2}.
\end{align*}
Using the above inequality  and the definition of $\mu$ of \eqref{nonlinear_map_mu}, we obtain
\begin{align}
	&C_1\|I_h^ku-\mu(\underline{\theta}_h)\|_{1,h}
	\leq\cN_h^{\rm lin}(u;I_h^ku-\underline{\theta}_h, \underline{v}_h)+ \cN_h(\underline{\theta}_h;\underline{v}_h)+ \tilde{C}_2\|I_h^ku-\mu(\underline{\theta}_h)\|_{L^2}.\label{garding_est_1}
\end{align}
Rewriting the first and second terms of the above equation, we obtain
\begin{align}
	\cN_h^{\rm lin}(u;I_h^ku-\underline{\theta}_h, \underline{v}_h)+ \cN_h(\underline{\theta}_h;\underline{v}_h)
	&=\left(\cN_h^{\rm lin}(u;I_h^ku-\underline{\theta}_h, \underline{v}_h)-\tilde{\cN}_h^{\rm lin}(I_h^ku;I_h^ku-\underline{\theta}_h, \underline{v}_h)\right)\notag\\
	&\qquad+\left(\tilde{\cN}_h^{\rm lin}(I_h^ku;I_h^ku-\underline{\theta}_h, \underline{v}_h)+\cN_h(\underline{\theta}_h;\underline{v}_h)\right).\label{nonLin_exp}
\end{align}
Now, we compute some residuals related to the nonlinear PDE \eqref{StrongNonlin_pde}. Multiplying $v_h$ and applying the integration by parts on \eqref{StrongNonlin_pde}, we have
\begin{align}
	0&=-\integ\nabla{\cdot}\ba(u,\nabla u)v_h\dx+\integ f(u,\nabla u)v_h\dx=-\sit\nabla{\cdot}\ba(u,\nabla u)v_h\dx+\integ f(u,\nabla u)v_h\dx\notag\\
	&=\sum_{T\in\cT_h}\left(\int_T \ba(u,\nabla u){\cdot}\nabla v_T\dx+\sum_{F\in\cF_T}\int_F( v_F- v_T)\ba(u,\nabla u){\cdot}\bfnTF\ds\right)+\integ f(u,\nabla u)v_h\dx.\label{a_res1}
\end{align}
The first two terms of the above equations are rewritten by some adjustment of terms and using the definition of gradient reconstructed operator \eqref{Grad_recons} as
\begin{align}
	&\sum_{T\in\cT_h}\left(\int_T \ba(u,\nabla u){\cdot}\nabla v_T\dx+\sum_{F\in\cF_T}\int_F( v_F- v_T)\ba(u,\nabla u){\cdot}\bfnTF\ds\right)\notag\\
	&=\sum_{T\in\cT_h}\int_T \ba(u,\nabla u){\cdot}\Gtk \underline{v}_T\dx +\sum_{F\in\cF_T}\int_F( v_F- v_T)(\ba(u,\nabla u)-\pi_T^k\ba(u,\nabla u)){\cdot}\bfnTF\ds\notag\\
	&=\sit \ba(\pi_T^k u,\Gtk I_T^k u){\cdot}\Gtk\underline{v}_T\dx+\sum_{T\in\cT_h}\int_T  \left(\ba(u,\nabla u)-\ba(\pi_T^k u,\Gtk I_T^k u)\right){\cdot}\Gtk\underline{v}_T\dx\notag\\
	&\qquad+\sum_{F\in\cF_T}\int_F( v_F- v_T)\left(\ba(u,\nabla u)-\pi_T^k\ba(u,\nabla u)\right){\cdot}\bfnTF\ds.\label{a_res2}
\end{align}
Combining the above two equations \eqref{a_res1}--\eqref{a_res2}, we obtain
\begin{align}
	0&=\cN_h(I_h^k u;\underline{v}_h)-s_h(I_h^k u,\underline{v}_h)+\sum_{T\in\cT_h}\int_T  \left(\ba(u,\nabla u)-\ba(\pi_T^k u,\Gtk I_T^k u)\right){\cdot}\Gtk\underline{v}_T\dx\notag\\
	&\qquad+\sum_{F\in\cF_T}\int_F( v_F- v_T)\left(\ba(u,\nabla u)-\pi_T^k\ba(u,\nabla u)\right){\cdot}\bfnTF\ds+ \sit \left(f(u,\nabla u)-f(\pi_T^ku,\Gtk I_T^ku)\right)v_h\dx.\label{load_integ}
\end{align}
Estimating all but the first term $\cN_h(I_h^k u;\underline{v}_h)$ using Lemma~\ref{lem:1st_remainder_est} and the estimate $s_h(I_h^k u,\underline{v}_h)\leq Ch^{r+1}\|\underline{v}_h\|_{1,h}$, we have
\begin{align}\label{res_Nh}
	\cN_h(I_h^k u;\underline{v}_h)\leq Ch^{r+1}\|\underline{v}_h\|_{1,h}.
\end{align}
Using Lemma~\ref{lem:diff_Lin}, Lemma~\ref{lem:lin_est},  Assumption N.5 and the estimate \eqref{res_Nh}, we obtain from \eqref{nonLin_exp} that
\begin{align}
	\cN_h^{\rm lin}(u;I_h^ku-\underline{\theta}_h, \underline{v}_h)+ \cN_h(\underline{\theta}_h;\underline{v}_h)&\leq C_{\ba,f} Ch^{-d/2}\|\underline{\theta}_h-I_h^k u\|_{1,h}^2\|\underline{v}_h\|_{1,h}\notag\\
	&\quad+Ch^{r+1}\|\underline{v}_h\|_{1,h} +Ch^{r+1-d/2}\|\underline{\theta}_h-I_h^k u\|_{1,h}\|\underline{v}_h\|_{1,h}.\label{mu_rhs}
\end{align}
Combining \eqref{garding_est_1} and \eqref{mu_rhs}, we have
\begin{align}
	C_1\|I_h^ku-\mu(\underline{\theta}_h)\|_{1,h}&\leq C_{\ba,f} Ch^{-d/2}\|\underline{\theta}_h-I_h^k u\|_{1,h}^2\|\underline{v}_h\|_{1,h}+Ch^{r+1}\|\underline{v}_h\|_{1,h}\notag\\
	&\qquad +Ch^{r+1-d/2}\|\underline{\theta}_h-I_h^k u\|_{1,h}\|\underline{v}_h\|_{1,h}+\tilde{C}_2\|I_h^ku-\mu(\underline{\theta}_h)\|_{L^2}.\label{gard_pre_final1}
\end{align}
Now, we estimate $\|I_h^ku-\mu(\underline{\theta}_h)\|_{L^2}$ using the following dual problem: given $\underline{q}_h= I_h^ku-\mu(\underline{\theta}_h)$, find $\underline{\phi}_h\in\underline{U}_{h,0}^k$ such that
\begin{equation}\label{dual_l2_est}
	\cN_h^{\rm lin}(u; \underline{v}_h,\underline{\phi}_h)= (q_h, v_h)\quad \forall\underline{v}_h\in \underline{U}_{h,0}^k.
\end{equation}
Choosing $\underline{v}_h= I_h^ku-\mu(\underline{\theta}_h)$ in the above equation, using the definition \eqref{nonlinear_map_mu} and the estimate \eqref{mu_rhs}, we obtain
\begin{align*}
	&\|I_h^ku-\mu(\underline{\theta}_h)\|_{L^2}^2=\cN_h^{\rm lin}(u; I_h^ku-\mu(\underline{\theta}_h),\underline{\phi}_h)=\cN_h^{\rm lin}(u;I_h^ku-\underline{\theta}_h, \underline{\phi}_h)+ \cN_h(\underline{\theta}_h; \underline{\phi}_h)\\
	&\leq C_{\ba,f} Ch^{-d/2}\|\underline{\theta}_h-I_h^k u\|_{1,h}^2\|\underline{\phi}_h\|_{1,h}+Ch^{r+1}\|\underline{\phi}_h\|_{1,h}+Ch^{r+1-d/2}\|\underline{\theta}_h-I_h^k u\|_{1,h}\|\underline{\phi}_h\|_{1,h}.
\end{align*}
Using the a priori bound $\|\underline{\phi}_h\|_{1,h}\leq C\|I_h^ku-\mu(\underline{\theta}_h)\|_{L^2}$ of \eqref{dual_l2_est} (see  \eqref{dis_apriori}), we obtain
\begin{align}
	\|I_h^ku-\mu(\underline{\theta}_h)\|_{L^2}\leq  C_{\ba,f} Ch^{-d/2}\|\underline{\theta}_h-I_h^k u\|_{1,h}^2+Ch^{r+1} +Ch^{r+1-d/2}\|\underline{\theta}_h-I_h^k u\|_{1,h}.\label{gard_pre_final2}
\end{align}
Finally, use \eqref{gard_pre_final2} in \eqref{gard_pre_final1} and $\|\underline{v}_h\|_{1,h}=1$ to obtain 
\begin{align}
	\|I_h^ku-\mu(\underline{\theta}_h)\|_{1,h}\leq \tilde{C}\left(h^{r+1}+h^{r+1-d/2}\|\underline{\theta}_h-I_h^k u\|_{1,h}+h^{-d/2}\|\underline{\theta}_h-I_h^k u\|_{1,h}^2\right)\label{gard_final_est}
\end{align}
for some positive constant $\tilde{C}$ independent of $h$.
Choose $h_{*}$ such that
\begin{align*}
	(1+2\tilde{C}h_{*}^{r+1-d/2}+4\tilde{C}^2h_{*}^{r+1-d/2})\leq 2.
\end{align*}
This implies $(1+2\tilde{C}h^{r+1-d/2}+4\tilde{C}^2h^{r+1-d/2})\leq 2$ whenever $h\leq h_{*}$. Thus if $\|I_h^ku-\underline{\theta}_h\|_{1,h}\leq R(h):=2\tilde{C}h^{r+1}$, then  \eqref{gard_final_est} yields
\begin{align*}
	&\|I_h^ku-\mu(\underline{\theta}_h)\|_{1,h}\leq \tilde{C}\left(h^{r+1}+2\tilde{C}h^{2r+2-d/2}+4\tilde{C}^2h^{2r+2-d/2}\right)\\
	&\,\leq \tilde{C}h^{r+1}\left(1+2\tilde{C}h^{r+1-d/2}+4\tilde{C}^2h^{r+1-d/2}\right)\leq \tilde{C}h^{r+1}\times 2= R(h).
\end{align*}
Thus, for a sufficiently small $h$ $(h\leq h_{*})$, there exists a ball $B(I_h^ku;R(h))$ of radius $R(h)=2\tilde{C}h^{r+1}$ with center at $I_h^ku$ such that the following result holds
$$\|I_h^ku-\underline{\theta}_h\|_{1,h}\leq R(h)\Rightarrow \|I_h^ku-\mu(\underline{\theta}_h)\|_{1,h}\leq R(h).$$
Therefore, $\mu$ is a map from a closed and bounded (compact) convex ball to itself. Therefore, using the Brouwer fixed point theorem, it has a fixed point. This completes the proof.
\end{proof}

\begin{remark}
It can be observed that the requirement of the regularity assumption $u\in H^3(\Omega)$ is merely to have $(1+2\tilde{C}h^{r+1-d/2}+4\tilde{C}^2h^{r+1-d/2})\leq 2$ for a sufficiently small meshsize $h$. This can also be done under the less regularity assumption $u\in H^{2+\epsilon}(\Omega)$ when $d=2$ and $u\in H^{5/2+\epsilon}(\Omega)$ when $d=3$, for any $\epsilon>0$  so that $(r+1-d/2)>0$ for real number $r=\epsilon$ if $d=2$ and $r=1/2+\epsilon$ if $d=3$.
\end{remark}
We show the contraction result to prove the unique fixed point of $\mu$. Recall Assumption N.1--N.5, then the contraction result holds:
\begin{theorem}[Contraction result]\label{thm_cont_mu}
Adopt the aforementioned Assumptions~N.1--N.5. Let $u\in H^1_0(\Omega)$ be a solution to \eqref{StrongNonlin_pde}. Assume $u\in H^{r+2}(\cT_h)$ and $\ba(x,y,\bz)$ to be $(r+1)$ times continuously differentiable with respect to $x$, for some $r\in\{1,\ldots,k\}$.
Let $\underline{\theta}_1,\, \underline{\theta}_2\in B(I_h^ku;R(h))$. For sufficiently small $h$, the following contraction result holds:
$$\|\mu(\underline{\theta}_1)-\mu(\underline{\theta}_2)\|_{1,h}\leq Ch^{r+1-d/2} \|\underline{\theta}_1-\underline{\theta}_2\|_{1,h}.$$
\end{theorem}

\begin{proof}
For $\underline{\theta}_1,\, \underline{\theta}_2\in B(I_h^ku;R(h))$, $\mu(\underline{\theta}_1)$ and $\mu(\underline{\theta}_2)$ satisfy \eqref{nonlinear_map_mu}. That is
\begin{align}
	\cN_h^{\rm lin}(u;I_h^ku-\mu(\underline{\theta}_1), \underline{v}_h)= \cN_h^{\rm lin}(u;I_h^ku-\underline{\theta}_1, \underline{v}_h)+ \cN_h(\underline{\theta}_1; \underline{v}_h)\label{defn_mu_1}\\
	\cN_h^{\rm lin}(u;I_h^ku-\mu(\underline{\theta}_2), \underline{v}_h)= \cN_h^{\rm lin}(u;I_h^ku-\underline{\theta}_2, \underline{v}_h)+ \cN_h(\underline{\theta}_2; \underline{v}_h).\label{defn_mu_2}
\end{align}
Choose $\underline{\theta}_h=\mu(\underline{\theta}_1)-\mu(\underline{\theta}_2)$ in the \GardingsType inequality \eqref{Garding_infsup} to obtain
\begin{align}
	C_1\|\mu(\underline{\theta}_1)-\mu(\underline{\theta}_2)\|_{1,h}\leq \cN_h^{\rm lin}(u;\mu(\underline{\theta}_1)-\mu(\underline{\theta}_2), \underline{v}_h)+ \tilde{C}_2\|\mu(\underline{\theta}_1)-\mu(\underline{\theta}_2)\|_{L^2},\label{mu_garding}
\end{align}
for some $\|\underline{v}_h\|_{1,h}=1$.
From the definition of $\mu$ and subtracting \eqref{defn_mu_1} with \eqref{defn_mu_2}, we get
\begin{align}
	\cN_h^{\rm lin}(u;\mu(\underline{\theta}_2)-\mu(\underline{\theta}_1), \underline{v}_h)= \cN_h^{\rm lin}(u;\underline{\theta}_2-\underline{\theta}_1, \underline{v}_h)+ \cN_h(\underline{\theta}_1;\underline{v}_h)-\cN_h(\underline{\theta}_2; \underline{v}_h).\label{mu_diff}
\end{align}
Using the definitions of $\cN_h$ and $\tilde{\cN}_h^{\rm lin}$ and the Taylor's formula \eqref{2nd_Order_Taylor}, the last two terms of the above equation \eqref{mu_diff} yield
\begin{align}
	&\cN_h(\underline{\theta}_1;\underline{v}_h)- \cN_h(\underline{\theta}_2;\underline{v}_h)=\left(\cN_h(\underline{\theta}_1;\underline{v}_h)-\cN_h(I_h^k u;\underline{v}_h)\right)-\left( \cN_h(\underline{\theta}_2;\underline{v}_h)-\cN_h(I_h^k u;\underline{v}_h)\right)\notag\\
	&=\tilde{\cN}_h^{\rm lin}(I_h^ku;\underline{\theta}_1-\underline{\theta}_2, \underline{v}_h)+\sit \tilde{R}_{\ba}(\theta_{1T}-\pi_T^k u, \Gtk(\underline{\theta}_{1T}-I_T^k u)){\cdot}\Gtk\underline{v}_T\dx\notag\\
	&\quad+\sit \tilde{R}_{f}(\theta_{1T}-\pi_T^k u,\Gtk(\underline{\theta}_{1T}-I_T^k u))\,v_T\dx-\sit \tilde{R}_{\ba}(\theta_{2T}-\pi_T^k u,\Gtk(\underline{\theta}_{2T}-I_T^k u)){\cdot}\Gtk\underline{v}_T\dx\notag\\
	&\quad\qquad-\sit \tilde{R}_{f}(\theta_{2T}-\pi_T^k u,\Gtk(\underline{\theta}_{2T}-I_T^k u))\,v_T\dx.\label{Diff_Nh_theta}
\end{align}
To obtain a difference term of the form $(\underline{\theta}_2-\underline{\theta}_1)$ from the last four terms of the above expression~\eqref{Diff_Nh_theta}, we use the definition of the residuals $\tilde{R}_{\ba}$ and $\tilde{R}_{f}$. Set $\underline{\xi}_1:=\underline{\theta}_1-I_h^k u, \underline{\xi}_2:=\underline{\theta}_2-I_h^k u$ and $\underline{\eta}:=\underline{\theta}_2-\underline{\theta}_1$. From the definition of residual in \eqref{Res_Ra}, we have 
\begin{align}
	&\tilde{R}_{\ba}(\xi_1,\Gtk \underline{\xi}_1)-\tilde{R}_{\ba}(\xi_2,\Gtk \underline{\xi}_2)\notag\\
	&=\left(\ba(\theta_1,\Gtk\underline{\theta}_1)-\ba(\pi_T^ku,\Gtk I_T^ku)+\ba_y(\pi_T^ku,\Gtk I_T^k u)(\pi_T^ku-\theta_{1T})+\ba_{\bz}(\pi_T^ku,\Gtk I_T^k u)\Gtk(I_T^ku-\underline{\theta}_{1T})\right) \notag\\
	&\quad-\left(\ba(\theta_2,\Gtk\underline{\theta}_2)-\ba(\pi_T^ku,\Gtk I_T^ku)+\ba_y(\pi_T^ku,\Gtk I_T^k u)(\pi_T^ku-\theta_{2T})+\ba_{\bz}(\pi_T^ku,\Gtk I_T^k u)\Gtk(I_T^ku-\underline{\theta}_{2T})\right)\notag\\
	&=\ba(\theta_1,\Gtk\underline{\theta}_1)-\ba(\theta_2,\Gtk\underline{\theta}_2)+\ba_y(\pi_T^ku,\Gtk I_T^k u)\eta_T+\ba_{\bz}(\pi_T^ku,\Gtk I_T^k u)\Gtk\underline{\eta}_{T}\notag\\
	&=\ba(\theta_1,\Gtk\underline{\theta}_1)-\ba(\theta_2,\Gtk\underline{\theta}_2)+\ba_y(\theta_{2T},\Gtk\underline{\theta}_{2T})\eta_T+\ba_{\bz}(\theta_{2T},\Gtk\underline{\theta}_{2T})\Gtk\underline{\eta}_{T}\notag\\
	&\quad+\left(\ba_y(\pi_T^ku,\Gtk I_T^k u)-\ba_y(\theta_{2T},\Gtk\underline{\theta}_{2T})\right)\eta_T+\left(\ba_{\bz}(\pi_T^ku,\Gtk I_T^k u)-\ba_{\bz}(\theta_{2T},\Gtk\underline{\theta}_{2T})\right)\Gtk\underline{\eta}_{T}\notag\\
	&=\tilde{R}_{\ba}(\eta_T,\Gtk\underline{\eta}_T)+\left(\ba_y(\pi_T^ku,\Gtk I_T^k u)-\ba_y(\theta_{2T},\Gtk\underline{\theta}_{2T})\right)\eta_T+\left(\ba_{\bz}(\pi_T^ku,\Gtk I_T^k u)-\ba_{\bz}(\theta_{2T},\Gtk\underline{\theta}_{2T})\right)\Gtk\underline{\eta}_{T}.\label{Rem_est_a}
\end{align}
Corollary~\ref{cor:2nd_remainder_est}, Assumption N.5, \eqref{Rem_est_a} and the triangle inequality $\|\underline{\theta}_2-\underline{\theta}_1\|_{1,h}\leq \|\underline{\theta}_2-I_h^k u\|_{1,h}+\|\underline{\theta}_1-I_h^k u\|_{1,h}$ lead to an estimate for the second and fourth terms of \eqref{Diff_Nh_theta} as 
\begin{align}
	&\sit \tilde{R}_{\ba}(\theta_{1T}-\pi_T^k u, \Gtk(\underline{\theta}_{1T}-I_T^k u)){\cdot}\Gtk\underline{v}_T\dx-\sit \tilde{R}_{\ba}(\theta_{2T}-\pi_T^k u,\Gtk(\underline{\theta}_{2T}-I_T^k u)){\cdot}\Gtk\underline{v}_T\dx\notag\\
	&\leq C C_{\ba} h^{-d/2}\|\underline{\theta}_{1}- \underline{\theta}_{2}\|_{1,h} \big{(}\|I_h^ku-\underline{\theta}_{1}\|_{1,h}+\|I_h^ku-\underline{\theta}_{2}\|_{1,h}\big{)}\|\underline{v}\|_{1,h}.
\end{align}
Exactly the same estimate holds for the combination of the third and fifth terms of \eqref{Diff_Nh_theta}.
Combining the above estimates, we obtain from \eqref{mu_diff} as
\begin{align}
	\cN_h^{\rm lin}(u;\mu(\underline{\theta}_2)-\mu(\underline{\theta}_1), \underline{v}_h)&\leq \left(\cN_h^{\rm lin}(u;\underline{\theta}_2-\underline{\theta}_1, \underline{v}_h)-\tilde{\cN}_h^{\rm lin}(I_h^ku;\underline{\theta}_{2}- \underline{\theta}_{1}, \underline{v}_h)\right)\notag\\
	&\quad+C C_{\ba,f} h^{-d/2}\|\underline{\theta}_{1}- \underline{\theta}_{2}\|_{1,h} \bigg{(}\|I_h^ku-\underline{\theta}_{1}\|_{1,h}+\|I_h^ku-\underline{\theta}_{2}\|_{1,h}\bigg{)}\|\underline{v}_h\|_{1,h}.\label{cts_lin_nonselfadj05}
\end{align}
Using Lemma~\ref{lem:diff_Lin}, we obtain from \eqref{cts_lin_nonselfadj05}
\begin{align}
	&\cN_h^{\rm lin}(u;\mu(\underline{\theta}_2)-\mu(\underline{\theta}_1), \underline{v}_h)\notag\\
	&\quad\leq Ch^{r+1-d/2}\|\underline{\theta}_{1}- \underline{\theta}_{2}\|_{1,h}\|\underline{v}_h\|_{1,h}+C C_{\ba,f} h^{-d/2}\|\underline{\theta}_{1}- \underline{\theta}_{2}\|_{1,h} \bigg{(}\|I_h^ku-\underline{\theta}_{1}\|_{1,h}+\|I_h^ku-\underline{\theta}_{2}\|_{1,h}\bigg{)}\|\underline{v}_h\|_{1,h}.\label{Nh_est_mu}
\end{align}
To obtain the estimate for $L^2$-term $\|\mu(\underline{\theta}_1)-\mu(\underline{\theta}_2)\|_{L^2}$, consider the dual linear problem: given $\underline{q}_h=\mu(\underline{\theta}_1)-\mu(\underline{\theta}_2)$, find $\underline{\phi}_h\in \underline{U}_{h,0}^k$ such that
\begin{align}
	\cN_h^{\rm lin}(u; \underline{v}_h,\underline{\phi}_h)= (q_h,v_h)\quad \forall \underline{v}_h\in \underline{U}_{h,0}^k.\label{dual_cont_prob} 
\end{align}
Choose $\underline{v}_h= \mu(\underline{\theta}_1)-\mu(\underline{\theta}_2)$ to obtain from \eqref{dual_cont_prob} and \eqref{Nh_est_mu}
\begin{align*}
	&\|\mu(\underline{\theta}_1)-\mu(\underline{\theta}_2)\|_{L^2}^2=\cN_h^{\rm lin}(u; \mu(\underline{\theta}_1)-\mu(\underline{\theta}_2),\underline{\phi}_h)\\
	&\leq C_{\ba,f} Ch^{-d/2}\|\underline{\theta}_{1}- \underline{\theta}_{2}\|_{1,h}\left(\|I_h^ku-\underline{\theta}_{1}\|_{1,h}+\|I_h^ku-\underline{\theta}_{2}\|_{1,h}\right)\|\underline{\phi}_h\|_{1,h}+Ch^{r+1-d/2}\|\underline{\theta}_{1}- \underline{\theta}_{2}\|_{1,h}\|\underline{\phi}_h\|_{1,h}.
\end{align*}
The a priori bound $\|\underline{\phi}_h\|_{1,h}\leq \|\mu(\underline{\theta}_1)-\mu(\underline{\theta}_2)\|_{L^2}$ of \eqref{dual_cont_prob} leads to
\begin{align}
	&\|\mu(\underline{\theta}_1)-\mu(\underline{\theta}_2)\|_{L^2}\leq C_{\ba,f} C\|\underline{\theta}_{1}- \underline{\theta}_{2}\|_{1,h}h^{-d/2}\left(\|I_h^ku-\underline{\theta}_{1}\|_{1,h}+\|I_h^ku-\underline{\theta}_{2}\|_{1,h}\right)+Ch^{r+1-d/2}\|\underline{\theta}_{1}- \underline{\theta}_{2}\|_{1,h}.\label{cts_lin_nonselfadj08}
\end{align}
Using the estimates \eqref{Nh_est_mu} with $\|\underline{v}_h\|_{1,h}=1$ and \eqref{cts_lin_nonselfadj08} in  \eqref{mu_garding}, we obtain
\begin{align*}
	&\|\mu(\underline{\theta}_1)-\mu(\underline{\theta}_2)\|_{1,h}\leq C_{\ba,f} C\|\underline{\theta}_{1}- \underline{\theta}_{2}\|_{1,h}h^{-d/2}\left(\|I_h^ku-\underline{\theta}_{1}\|_{1,h}+\|I_h^ku-\underline{\theta}_{2}\|_{1,h}\right)+Ch^{r+1-d/2}\|\underline{\theta}_{1}- \underline{\theta}_{2}\|_{1,h}.
\end{align*}
Since $\underline{\theta}_{1},\, \underline{\theta}_{2}\in B(I_h^ku;R(h))$ with $R(h)=2\tilde{C}h^{r+1}$, that is,
\begin{equation*}
	\|I_h^ku-\underline{\theta}_{1}\|_{1,h}\leq 2\tilde{C}h^{r+1}\quad \text{and}\quad \|I_h^ku-\underline{\theta}_{2}\|_{1,h}\leq 2\tilde{C}h^{r+1}.
\end{equation*}
For sufficiently small meshsize $h$, we have
\begin{equation*}
	\|\mu(\underline{\theta}_1)-\mu(\underline{\theta}_2)\|_{1,h}\leq Ch^{r+1-d/2}\|\underline{\theta}_{1}- \underline{\theta}_{2}\|_{1,h},    
\end{equation*}
for some positive constant independent of $h$.
This completes the proof.
\end{proof}

For sufficiently small $h$, the above Theorem~\ref{thm_cont_mu} proves the local uniqueness of the fixed point of $\mu$ and hence the local uniqueness of the solution to \eqref{hho_dis_quasi}.

Adding and subtracting $\Ghk I_h^k u$, using triangle inequality, the definition of norm $\|\bullet\|_{1,h}$ in \eqref{hho_gtk_norm} and Theorem~\ref{thm_fixed_point}, we have the following error estimate under Assumptions~N.1--N.5:
\begin{theorem}[Error estimate]\label{thm_err_est_strong}
Adopt the aforementioned Assumptions~N.1--N.5. Let $u\in H^1_0(\Omega)$ be the solution to nonlinear problem \eqref{StrongNonlin_pde} and $\underline{u}_h\in\underline{U}_{h,0}^k$ be the solution to the discrete problem \eqref{hho_dis_quasi}. Assume $u\in H^{r+2}(\cT_h)$ and $\ba(x,y,\bz)$ to be $(r+1)$-times continuously differentiable with respect to $x$, for some $r\in\{1,\ldots,k\}$. Then for sufficiently small $h$, we have
\begin{align}\label{err_est_quasi}
	\|\nabla u-\Ghk\underline{u}_h\|\leq Ch^{r+1},
\end{align}
for some positive constant $C$ independent of $h$.
\end{theorem}

\begin{remark}
	We observe that for the special case of the nonlinear function $\ba(x,u,\nabla)=a(x,u)\nabla u$, the authors in \cite{TG_GM_TP_22_HHO_Quasi} obtained optimal order error estimate for the lowest-order ($k=0$) HHO polynomial approximation. However, due to the strongly nonlinear problem, we obtain an optimal order error estimate for $k\geq 1$. The error estimate for the lowest-order polynomial approximations ($k=0$ for HHO and $k=1$ for various discontinuous Galerkin methods \cite{Gudi_NN_AKP_08_Strongly_nonlin,Bi_Cheng_Lin_16_Pointwise}) is still an open question that requires further study.
\end{remark}

\section{Numerical experiments}\label{sec:Numerics}
In this section, we perform some numerical experiments for the strongly nonlinear problem~\eqref{StrongNonlin_pde} using the HHO approximation described in \eqref{hho_dis_quasi}. Consider the following strongly nonlinear model problem \cite{Gudi_AKP_07_DG_quasi}:
\begin{subequations}\label{num_strong_nonlin}
\begin{align}
	-\nabla{\cdot}\left(\frac{\nabla u}{\sqrt{1+|\nabla u|^2}}\right) &= f\quad\text{in}~~ \Omega,\\
	u&=0\quad\text{on}~~ \partial\Omega,
\end{align}
\end{subequations}
where we have taken $\ba(x,y,\bz)=\bz(1+|\bz|^2)^{-1/2}$ and $f(x,y,\bz)=-f(x)$ in \eqref{StrongNonlin_pde} to obtain the above problem.
From the application point of view, the above model problem~\eqref{num_strong_nonlin} describes the mean curvature flow. We verify Assumptions~N.1--N.2 as follows: for $\bz=(z_1,z_2)$, we obtain the following derivative matrix
\begin{align}
	\ba_{\bz}(y,\bz)=\left[a^{ij}(x,y,\bz)\right]_{i,j=1}^{2}=\left[\frac{\partial a_i}{\partial z_j}\right]_{i,j=1}^{2}=R(\bz)\begin{bmatrix}
		1+z_2^2 & -z_1z_2\\
		-z_1z_2 & 1+z_1^2
	\end{bmatrix},
\end{align}
where $R(z)=(1+z_1^2+z_2^2)^{-3/2}$. The ellipticity condition \eqref{cts_ellipticity} of Assumption N.2 is verified as follows: for $\xi=(\xi_1,\xi_2)\in \mbR^2\setminus \boldsymbol{0}$,
\begin{align*}
	\sum_{i,j=1}^2 a^{ij}(y,z)\xi_i\xi_j=R(\bz)\left((1+z_2^2)\xi_1^2-2z_1z_2\xi_1\xi_2+(1+z_1^2)\xi_2^2\right)=R(\bz)\left((z_2\xi_1-z_1\xi_2)^2+(\xi_1^2+\xi_2^2)\right)
\end{align*}
Since $0\leq (z_2\xi_1-z_1\xi_2)^2\leq 2|\bz|^2|\xi|^2$, we have the following boundedness
\begin{align*}
	R(\bz)|\xi|^2\leq\sum_{i,j=1}^2 a^{ij}(y,\bz)\xi_i\xi_j\leq R(\bz)(1+2|\bz|^2)|\xi|^2.
\end{align*}
For numerical experiments, we consider the domain to be a unit square, i.e. $\Omega:=(0,1)\times (0,1)$. The source term $f$ is taken in such a way that the exact solution reads $u(x,y)=x(1-x)y(1-y)$. For $\bz=\nabla u$, $R(\bz)$ is bounded below by a constant $\lambda_0$. Since $R(\bz)\leq 1$, $R(\bz)(1+2|\bz|^2)$ is bounded above by a constant $\Lambda_0$. Since $u$ is sufficiently smooth, Assumptions~N.1--N.2 follow. Assumption~N.4 follows due to the smoothness of $u$. We can observe that $f_y(y,\bz)=f_{\bz}(y,\bz)=\nabla{\cdot}\ba_y(y,\bz)=0$. This verifies Assumption~N.4. In the numerical tests, we consider quasi-uniform mesh sequences that validate Assumption~N.5.

We describe an iterative step to obtain the discrete solution. The nonlinear map $\mu$ defined in \eqref{nonlinear_map_mu} helps to design an iterative process, where we replace the exact solution with the computed solution from the previous step. We start with an initial guess $u_h^0\in \underline{U}_{h,0}^k$ obtained from solving the Dirichlet Poisson problem $-\Delta u=f$ with the same load function $f$ as defined above. The $(n+1)$-th iteration is given by
\begin{align}
\tilde{\cN}_h^{\rm lin}(\underline{u}_h^{n};\underline{u}_h^{n+1},\underline{v}_h)=\tilde{\cN}_h^{\rm lin}(\underline{u}_h^{n};\underline{u}_h^{n},\underline{v}_h)-\cN_h(\underline{u}_h^n;\underline{v}_h)\fl \underline{v}_h \in\underline{U}_{h,0}^k,\, n=0,1,2,\ldots,
\end{align} 
where the linearized  $\tilde{\cN}_h^{\rm lin}$ and nonlinear $\cN_h$  forms are as defined in \eqref{defn_Bh_lin} and \eqref{hho_dis_quasi}, respectively. The stopping criterion is prescribed by a tolerance $10^{-8}$ for the difference of two successive iterative solutions as $\|\Ghk(\underline{u}_h^{n+1}-\underline{u}_h^{n})\|/\|\Ghk\underline{u}_h^{n+1}\|\leq 10^{-8}$.

\begin{figure}
\begin{center}
	\subfloat[]{\includegraphics[height=0.2\textwidth,width=0.25\textwidth]{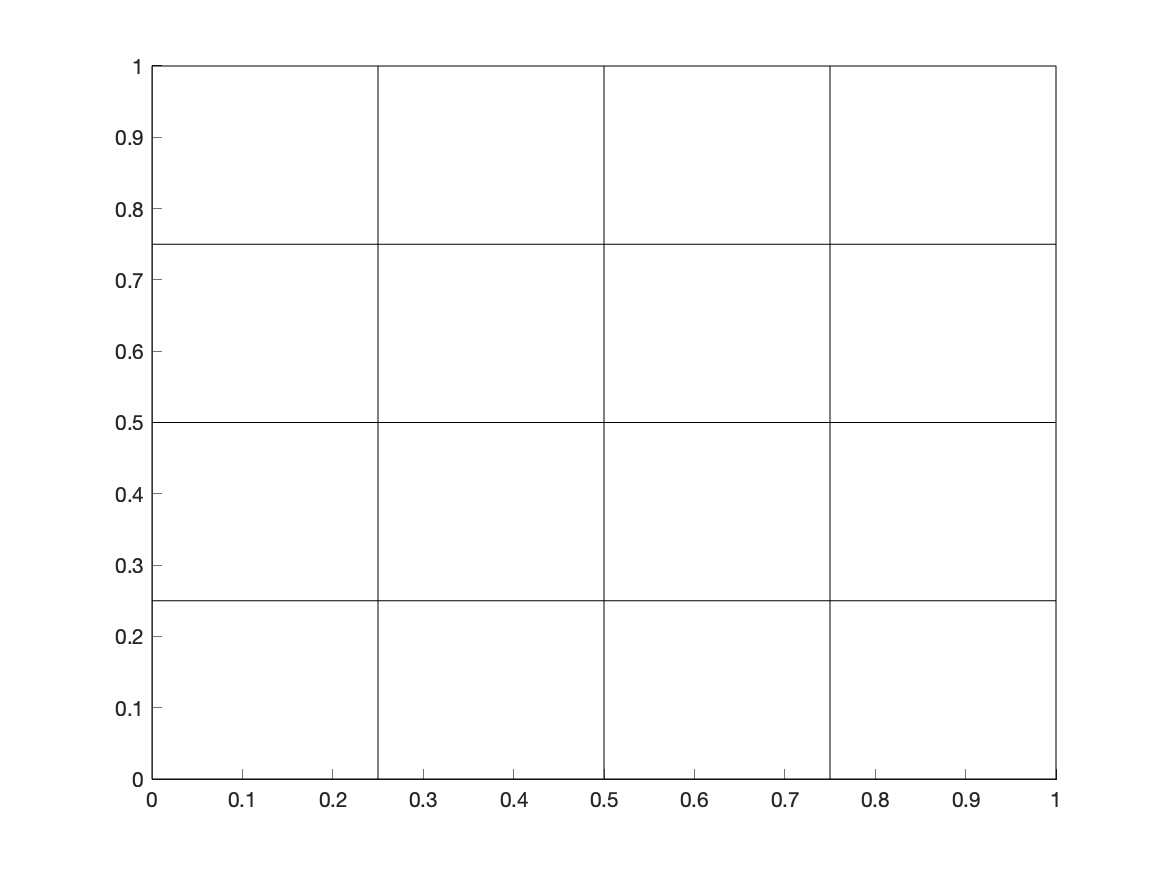}}
	\subfloat[]{\includegraphics[height=0.2\textwidth,width=0.25\textwidth]{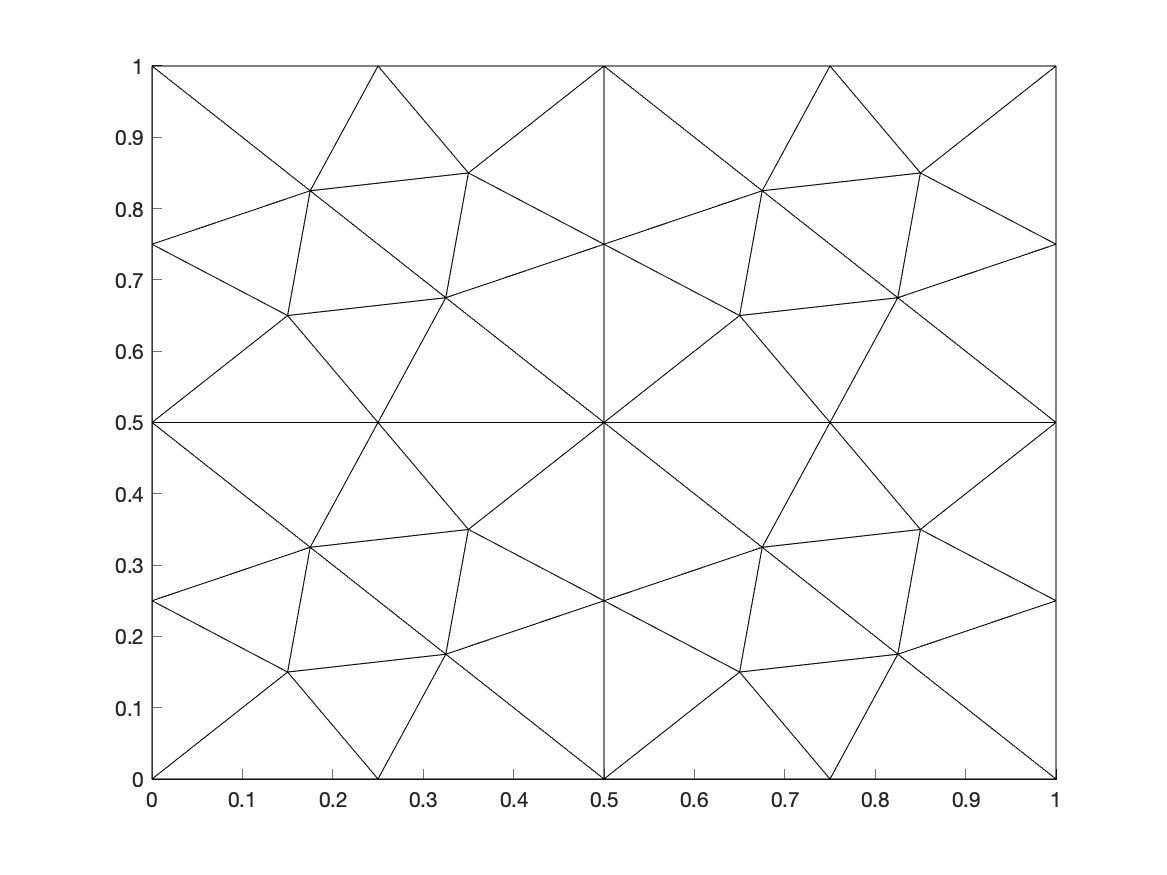}}
	\subfloat[]{\includegraphics[height=0.2\textwidth,width=0.25\textwidth]{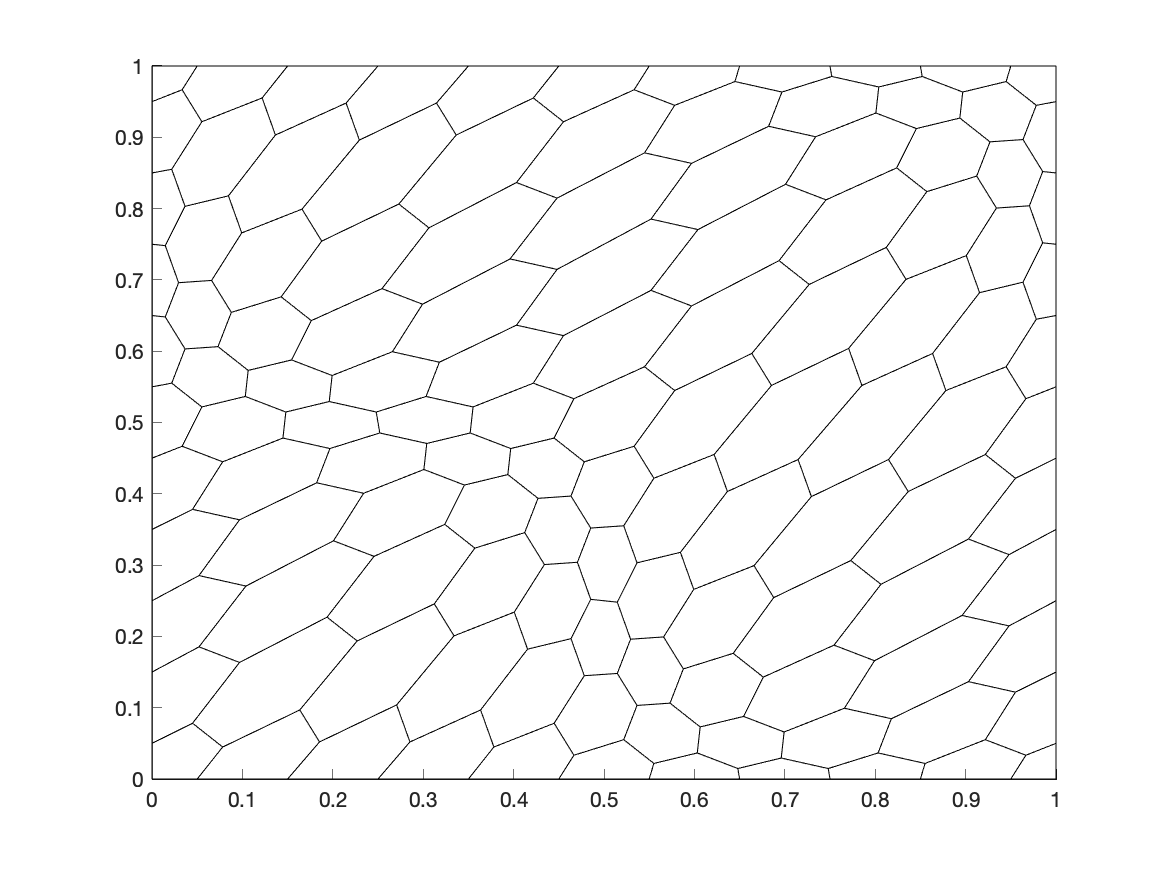}}
	\subfloat[]{\includegraphics[height=0.2\textwidth,width=0.25\textwidth]{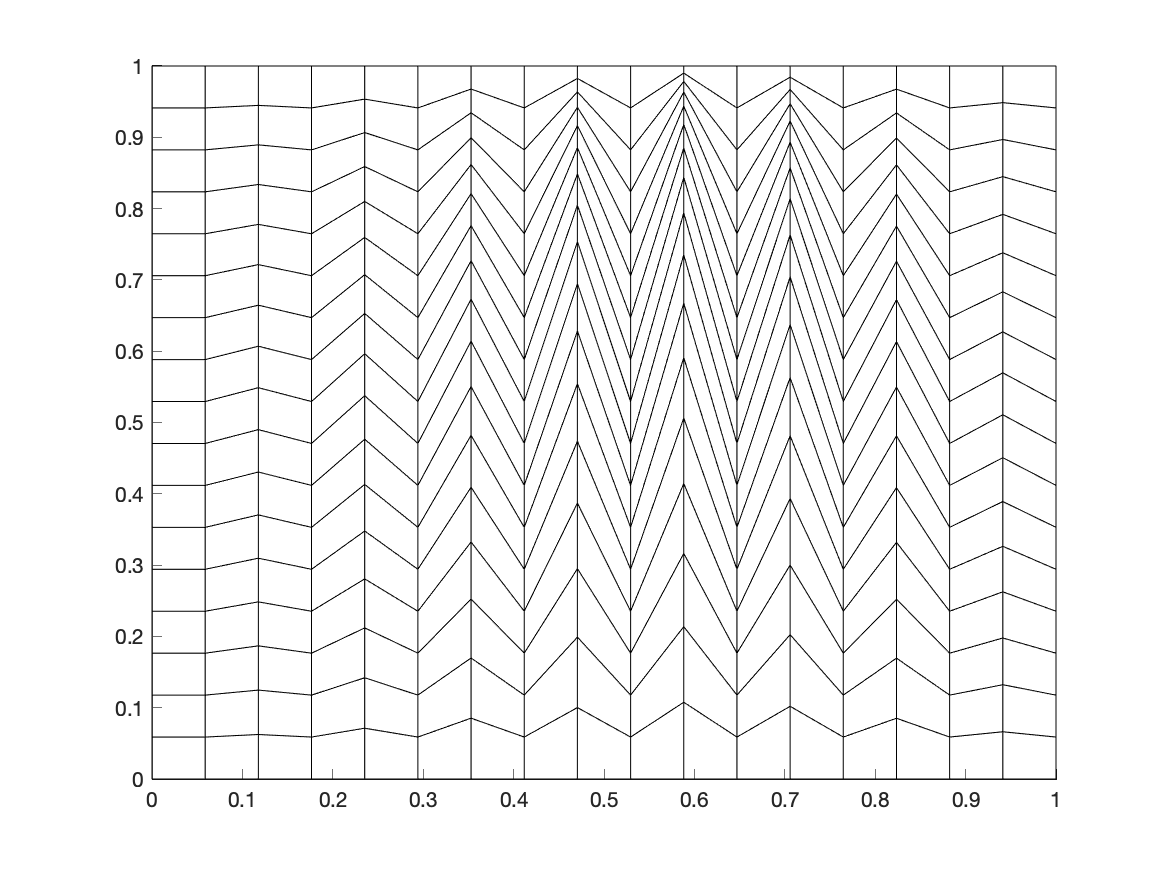}}
	\caption{(a) Cartesian, (b) Triangular, (c) hexagonal and (d) Kershaw initial meshes \cite{TG_GM_TP_22_HHO_Quasi}.}
	\label{fig:Initial_Meshes}
\end{center}
\end{figure}

We perform numerical tests on four different families of meshes: Cartesian, triangular, hexagonal and Kershaw meshes. Their initial meshes are shown in Figure~\ref{fig:Initial_Meshes}. For details on the mesh families, we refer \cite{Her_Hub_FVCA_mesh_08} to the Cartesian, triangular and Kershaw mesh families and \cite{Piet_Lem_HexaMesh_15} the hexagonal mesh family.
We adapt some of the basic implementation methodologies for the HHO methods from \cite{Piet_Jero_HHO_Book_20,Cicu_Piet_Ern_18_Comp_HHO,Piet_Ern_Lem_14_arb_local}.
It has been observed that the iterative step terminates within $4$ steps using the above stopping criterion.
The empirical rate of convergence is given by
\begin{align*}
\texttt{rate}(\ell):=\log \big(e_{h_{\ell}}/e_{h_{\ell-1}} \big)/\log \big(h_\ell/h_{\ell-1} \big)\text{ for } \ell=1,2,3,\ldots,
\end{align*}
where $e_{h_{\ell}}$ and $e_{h_{\ell-1}}$ are the errors associated to the two consecutive meshsizes $h_\ell$ and $h_{\ell-1}$, respectively.

\begin{figure}
\begin{center}
	\subfloat[]{\includegraphics[height=0.4\textwidth,width=0.5\textwidth]{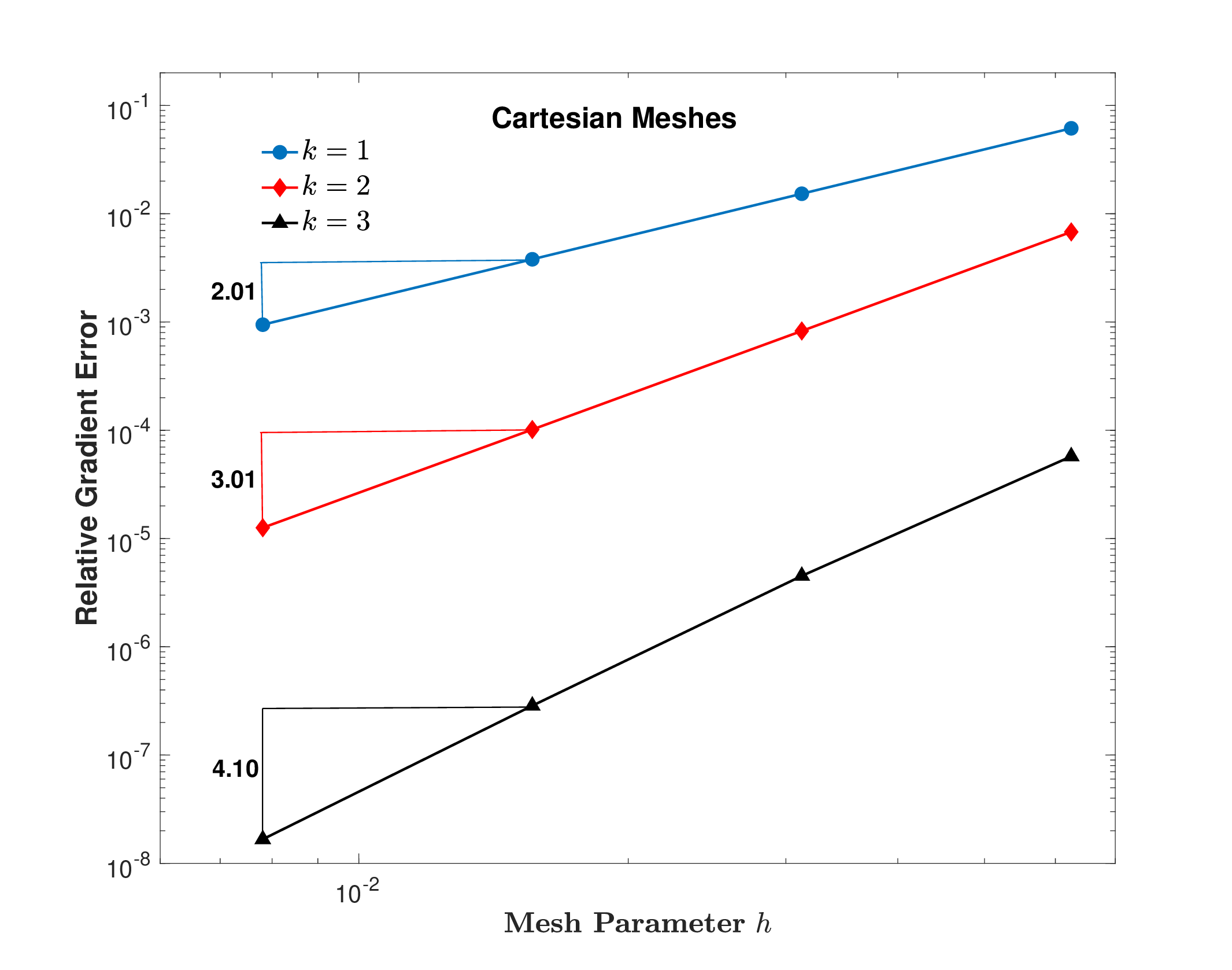}}
	\subfloat[]{\includegraphics[height=0.4\textwidth,width=0.5\textwidth]{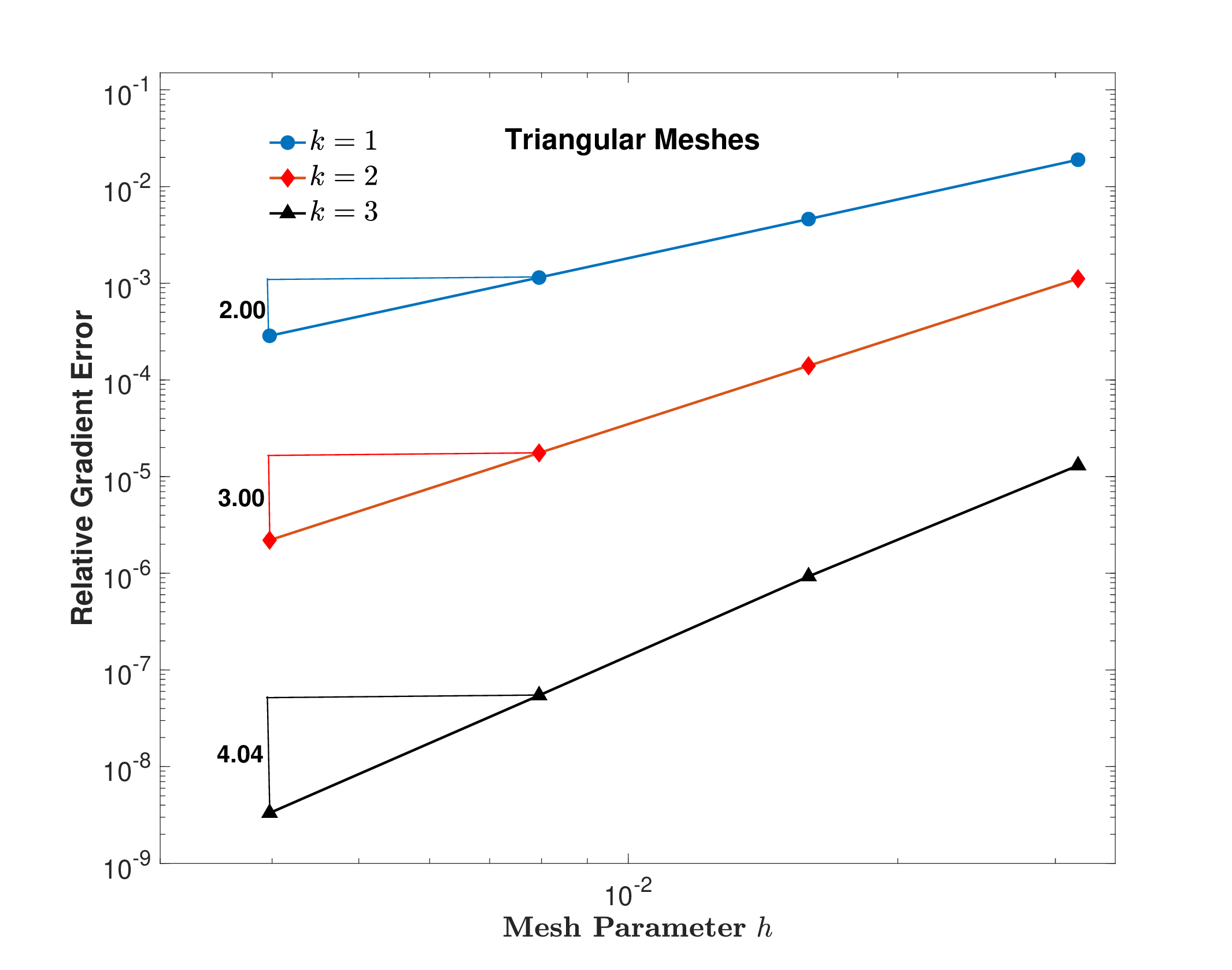}}\\
	\subfloat[]{\includegraphics[height=0.4\textwidth,width=0.5\textwidth]{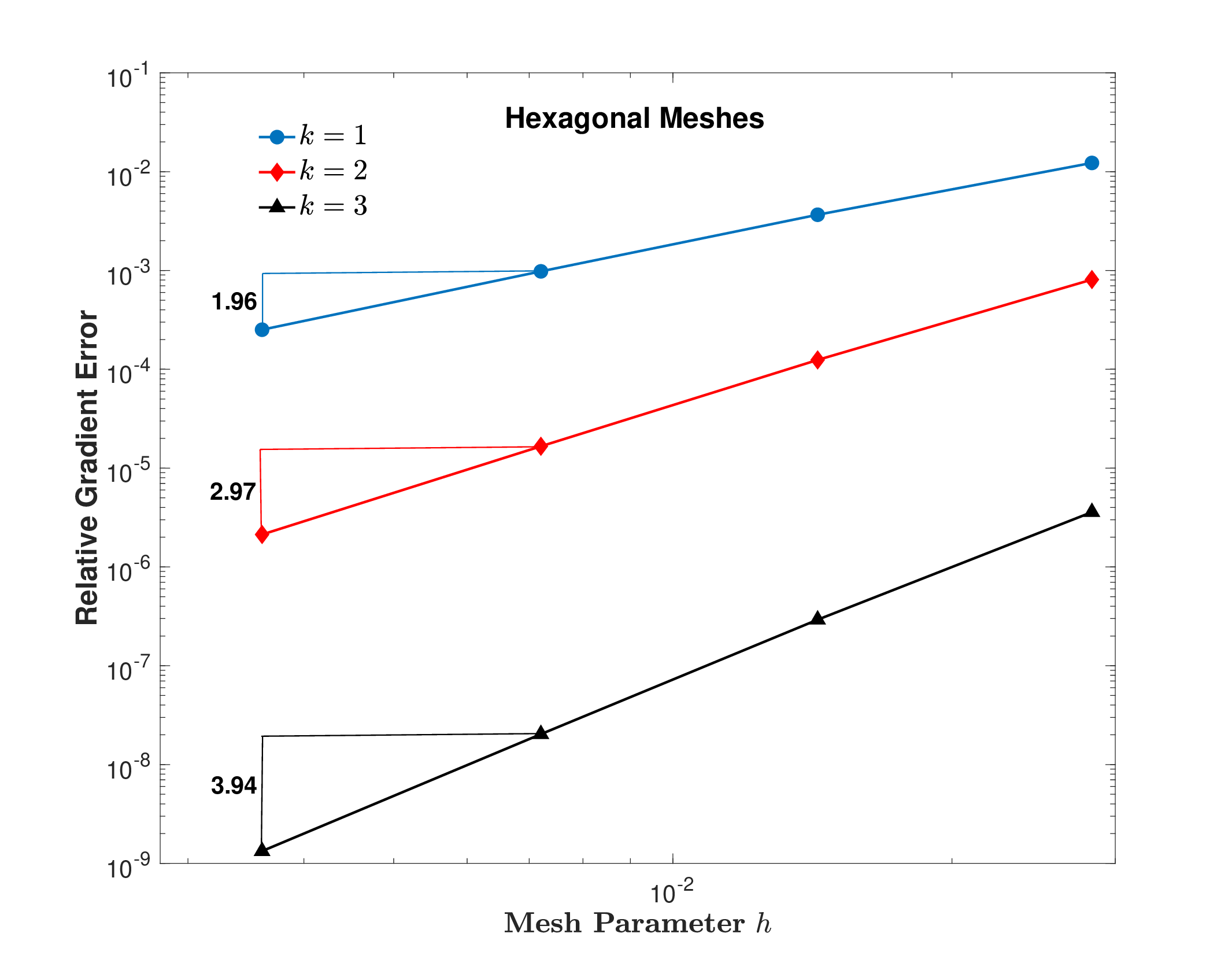}}
	\subfloat[]{\includegraphics[height=0.4\textwidth,width=0.5\textwidth]{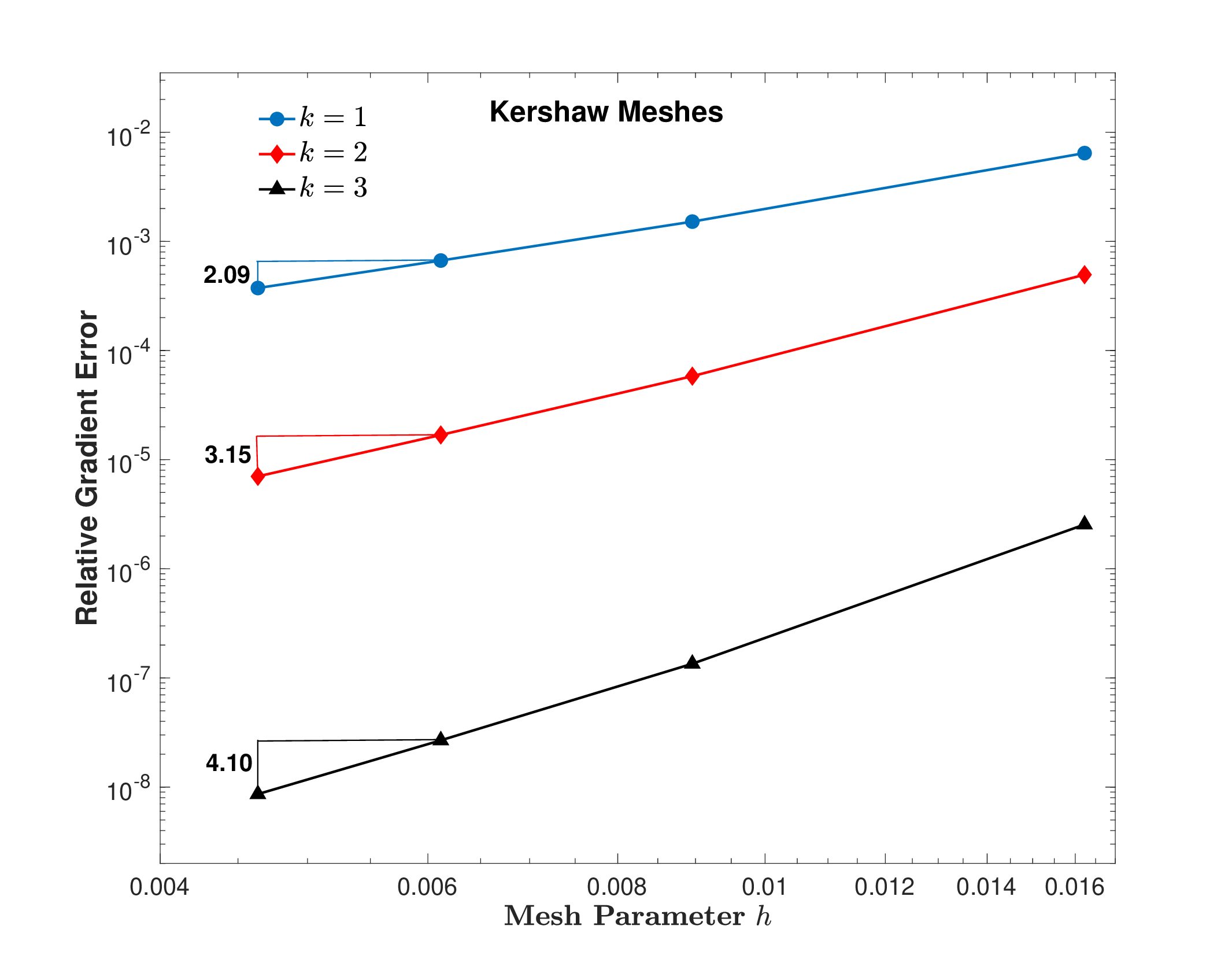}}
	\caption{Convergence histories for the relative gradient error on the (a) Cartesian, (b) Triangular,  (c) hexagonal and (d) Kershaw meshes.}
	\label{fig:Strong_Conv_His_Triag_Cart}
\end{center}
\end{figure}

In Table~\ref{table:Strong_Triag}--\ref{table:Strong_Kers}, we have shown the relative gradient error $e_h=\|\nabla u-\Ghk\underline{u}_h\|/\|\nabla u\|$ and its convergence rate for the Cartesian, triangular, hexagonal and Kershaw mesh families. The convergence histories for the relative gradient error $e_h$ with respect to meshsize $h$ have been plotted in Figure~\ref{fig:Strong_Conv_His_Triag_Cart}, where we have considered the Cartesian, triangular, hexagonal and Kershaw meshes for the polynomial degree $k=1,2,3$.
The empirical rates of convergence for the polynomial degree $k=1,2,3$ are close to $2,3,4$ for each mesh family. The empirical convergence rates obey the theoretical convergence rate of Theorem~\ref{thm_err_est_strong}.

\begin{table}
\caption{Relative gradient errors and convergence rates on the Cartesian meshes.}\label{table:Strong_Cart}
\begin{center}
	\begin{tabular}{  c c  c c c c c }     
		\hline
		\multirow{2}{*}{$h$}  & \multicolumn{2}{c}{$k=1$} &\multicolumn{2}{c}{$k=2$} &\multicolumn{2}{c}{$k=3$}\\
		\cline{2-7}
		&$e_h$ & \texttt{rate}  & $e_h$ & \texttt{rate} & $e_h$ & \texttt{rate}\\ 
		\hline 
		0.0625   & 0.6150e--1  & --     & 0.6791e--2 & --     & 0.5741e--4 & --\\  
		0.0313   & 0.1529e--1  & 2.008  & 0.8262e--3 & 3.039  & 0.4518e--5 & 3.668 \\ 
		0.0156   & 0.3795e--2  & 2.011  & 0.1015e--3 & 3.024  & 0.2857e--6 & 3.983 \\ 
		0.0078   & 0.9442e--3  & 2.007  & 0.1258e--4 & 3.013  & 0.1669e--7 & 4.098 \\ \hline
	\end{tabular}
\end{center}
\end{table}

\begin{table}
\caption{Relative gradient errors and convergence rates on the triangular meshes.}\label{table:Strong_Triag}
\begin{center}
	\begin{tabular}{  c c  c c c c c }     
		\hline
		\multirow{2}{*}{$h$}  & \multicolumn{2}{c}{$k=1$} &\multicolumn{2}{c}{$k=2$} &\multicolumn{2}{c}{$k=3$}\\
		\cline{2-7}
		&$e_h$ & \texttt{rate}  & $e_h$ & \texttt{rate} & $e_h$ & \texttt{rate}\\ 
		\hline 
		0.0318   & 0.1894e--1  & --     & 0.1113e--2 & --     & 0.1303e--4 & --\\  
		0.0159   & 0.4611e--2  & 2.039  & 0.1400e--3 & 2.991  & 0.9280e--6 & 3.812 \\ 
		0.0080   & 0.1145e--2  & 2.009  & 0.1756e--4 & 2.994  & 0.5459e--7 & 4.087 \\ 
		0.0040   & 0.2860e--3  & 2.002  & 0.2199e--5 & 2.998  & 0.3316e--8 & 4.041 \\ \hline
	\end{tabular}
\end{center}
\end{table}

\begin{table}
\caption{Relative gradient errors and convergence rates on the hexagonal meshes.}\label{table:Strong_Hexa}
\begin{center}
	\begin{tabular}{  c c  c c c c c }     
		\hline
		\multirow{2}{*}{$h$}  & \multicolumn{2}{c}{$k=1$} &\multicolumn{2}{c}{$k=2$} &\multicolumn{2}{c}{$k=3$}\\
		\cline{2-7}
		&$e_h$ & \texttt{rate}  & $e_h$ & \texttt{rate} & $e_h$ & \texttt{rate}\\ 
		\hline 
		0.0283   & 0.1226e--1  & --     & 0.8093e--3 & --     & 0.3585e--5 & --\\  
		0.0143   & 0.3665e--2  & 1.773  & 0.1243e--3 & 2.750  & 0.2923e--6 & 3.680 \\ 
		0.0072   & 0.9796e--3  & 1.920  & 0.1663e--4 & 2.928  & 0.2030e--7 & 3.882 \\ 
		0.0036   & 0.2515e--3  & 1.964  & 0.2127e--5 & 2.970  & 0.1330e--8 & 3.937 \\ \hline
	\end{tabular}
\end{center}
\end{table}

\begin{table}
\caption{Relative gradient errors and convergence rates on the Kershaw meshes.}\label{table:Strong_Kers}
\begin{center}
	\begin{tabular}{  c c  c c c c c }     
		\hline
		\multirow{2}{*}{$h$}  & \multicolumn{2}{c}{$k=1$} &\multicolumn{2}{c}{$k=2$} &\multicolumn{2}{c}{$k=3$}\\
		\cline{2-7}
		&$e_h$ & \texttt{rate}  & $e_h$ & \texttt{rate} & $e_h$ & \texttt{rate}\\ 
		\hline 
		0.0162   & 0.6439e--2  & --     & 0.4950e--3 & --     & 0.2544e--5 & --\\  
		0.0080   & 0.1517e--2  & 2.433  & 0.5818e--4 & 3.603  & 0.1349e--6 & 4.943 \\ 
		0.0061   & 0.6672e--3  & 2.156  & 0.1684e--4 & 3.255  & 0.2681e--7 & 4.240 \\ 
		0.0046   & 0.3737e--3  & 2.091  & 0.7026e--5 & 3.154  & 0.8617e--8 & 4.095 \\ \hline
	\end{tabular}
\end{center}
\end{table}

\section{Conclusions}\label{sec:Conclusion}
In this article, we studied the HHO finite element approximation for a class of strongly nonlinear elliptic PDEs. We proved the well-posedness of a discrete linearized problem using the \Gardings type inequality, where the lower-order $L^2$-term has been controlled by some estimates of the continuous linearized problem. We adapted the methodology of the fixed point arguments and the contraction principle in order to establish the existence of a discrete local solution. We obtained the optimal order error estimate in the energy norm as a by-product of the analysis. Several numerical experiments are performed to illustrate the optimal rate of convergence.
	
\bibliographystyle{siam}
\bibliography{HHO_Bib}
	
\end{document}